\documentclass[12pt,a4paper,leqno]{amsart}
\usepackage{amsmath}
\usepackage{amssymb}

\newtheorem{theorem}{Theorem}[section]

\newtheorem*{theoremC}{Theorem  C} 
 
\newtheorem{lemma}[theorem]{Lemma}
\newtheorem{proposition}[theorem]{Proposition}
\newtheorem{corollary}[theorem]{Corollary}

\newcommand{\C}{{\bf C}}       
\newcommand{\B}{{\bf B}}       
\newcommand{\R}{{\bf R}}       

\catcode`\@=11 \@addtoreset{equation}{section} \catcode`\@=12

\title[Carleson measures for $H_s^p(w)$] { Carleson measures  for weighted  Hardy-Sobolev  spaces }
\author[Cascante]{Carme Cascante}
\address{Departament de Matem\`atica Aplicada i An\`alisi,
Facultat de Matem\`atiques, Universitat de Barcelona, Gran Via
585, 08071~Barcelona, Spain} \email{cascante@ub.edu}

\author[Ortega]{Joaquin M. Ortega }
\address{ Departament de Matem\`atica Aplicada i An\`alisi,
Facultat de Matem\`atiques,  Universitat de Barcelona, Gran Via
585, 08071~Barcelona, Spain} \email{ortega@ub.edu}

\thanks{Both authors partially supported by DGICYT Grant
MTM2005-08984-C02-02, and DURSI Grant 2005SGR 00611
.}

\date{}

\begin{document}
\begin{abstract}
We obtain characterizations of positive Borel measures $\mu$ on
$\B^n$
 so that some weighted  Hardy-Sobolev  are imbedded in  
  $L^p(d\mu)$,  where $w$ is an $A_p$ weight in the unit sphere  of $\C^n$. 
   
\end{abstract}
\maketitle

\section{Introduction}  \label{introduction}

 The  purpose of  this paper is the study of the positive Borel measures $\mu$ on ${\bf S}^n$, the   unit sphere in $\C^n$,   for which the weighted Hardy-Sobolev space $H_s^p(w)$   is imbedded in $L^p(d\mu)$, that is,   the Carleson measures for $H_s^p(w)$.

 The weighted Hardy-Sobolev space $H_s^p(w)$,   $0<s, p<+\infty$, consists of those functions $f$ holomorphic in $\B^n$ such that if ${\displaystyle f(z)=\sum_k f_k(z)}$ is its homogeneous polynomial expansion, and   ${\displaystyle (I+R)^s f(z)=\sum_k (1+k)^s f_k(z)},$  we have that
$${\displaystyle ||f||_{H_s^p(w)}=\sup_{0<r<1} ||(I+R)^s f_r||_{L^p(w)}<+\infty},$$ where $f_r(\zeta)=f(r\zeta)$.

We will consider weights $w$ in $A_p$ classes in ${\bf S}^n$, $1<p<+\infty$,  that is, weights in  ${\bf S}^n$ satisfying that there exists $C>0$ such that for any nonisotropic ball $B\subset {\bf S}^n$, $B=B(\zeta,r)=\{ \eta\in{\bf S}^n\,;\, |1-\zeta\overline{\eta}|<r\,\}$,
$$\left( \frac{1}{|B|}\int_B w d\sigma\right)\left( \frac{1}{|B|}\int_B w^{\frac{-1}{p-1}}d\sigma \right)^{p-1}\leq C,$$
where $\sigma$ is the Lebesgue measure on ${\bf S}^n$ and $|B|$ the Lebesgue measure of $B$. We will use the notation $\zeta\overline{\eta}$ to indicate the complex inner product in $\C^n$ given by ${\displaystyle{\zeta\overline{\eta}= \sum_{i=1}^n \zeta_i\overline{\eta_i}}}$, if $\zeta=(\zeta_1,\dots,\zeta_n)$, $\eta=(\eta_1,\dots,\eta_n)$.

  If $0<s<n$, any function $f$ in $H_s^p(w)$ can be expressed as $${\displaystyle{f(z)=C_s(g)(z):= \int_{{\bf S}^n} \frac{g(\zeta)}{(1-z\overline{\zeta})^{n-s}}d\sigma(\zeta)}},$$    where $d\sigma$ is the normalized Lebesgue measure on the unit sphere ${\bf S}^n$ and $g\in L^p(w)$, and    consequently,    $\mu$ is Carleson for $H_s^p(w)$ if there exists $C>0$ such that
 $$||C_s f||_{L^p(d\mu)}\leq C||f||_{L^p(w)}.$$

We denote by $K_s$ the nonisotropic potential operator defined by
$$K_s[f](z)=\int_{{\bf S}^n} \frac{f(\eta)}{|1-z\overline{\eta}|^{n-s}}d\sigma(\eta), \,\, z\in\overline{\B}^n.$$

The problem of characterizing the positive Borel measures $\mu$ on $\B^n$ for which there exists $C>0$ such that 
\begin{equation} \label{carlesonwithoutweights}
||K_s[f]||_{L^p(d\mu)}\leq C||f||_{L^p(d\sigma)},
\end{equation}
that is, the characterization of the Carleson measures for the space $K_s[L^p(d\sigma)]$ has been very well studied  and there exist different characterizations (see  for instance \cite{mazya}, \cite{adamshedberg} \cite{kermansawyer}).

The representation of the functions in $H_s^p$ in terms of the operator $C_s$ gives  that in dimension $1$ the Carleson measures for $K_s[L^p(d\sigma)]$ coincide with the Carleson measures for the Hardy-Sobolev space $H_s^p$ simply because the real part of $\frac1{(1-z\overline{\zeta})^{1-s}}$ is equivalent to $\frac1{|1-z\overline{\zeta}|^{1-s}}$.  This representation also shows   that in any dimension  every Carleson measure for $K_s[L^p(d\sigma)]$ is also a Carleson measure for $H_s^p$. The coincidence fails to be true  for $n>1$ in general, as it is shown in \cite{aherncohn} ( see also \cite{cascanteortega2}).

   Of course, when $n-sp<0$, the space $H_s^p$ consists of continuous functions on $\overline{\B}^n$, and in particular, the Carleson measures in this case are just the finite measures. But for $n-sp\geq 0$, and $n>1$, the characterization of the Carleson measures for $H_s^p$  still remains open. In the case where we are  ''near'' the regular case, that is when $n-sp<1$    it is shown in \cite{aherncohn}, \cite{cohnverbitsky2} and \cite{cohnverbitsky1}, that the Carleson measures for $H_s^p$ and $K_s[L^p(d\sigma)]$ are the same, and any of the different characterizations of the Carleson measures for the last ones also hold for $H_s^p$. 

 One of the main purposes of this paper is to extend this situation to  $H_s^p(w)$ for $w$ a weight in   $A_p$.  If $E\subset{\bf S}^n$ is measurable, we define
 $$W(E)=\int_E wd\sigma.$$ 
 A weight $w$   satisfies a doubling condition of order $\tau$,  if  there exists $\tau>0$ such that for any nonisotropic ball $B$ in ${\bf S}^n$, $W(2^kB)\leq C 2^{k\tau}W(B)$.
 
  It is well known that any weight in $A_p$ satisfies a doubling condition of some order $\tau$ strictly less than $np$. We begin   observing that if $\tau-sp<0$, the space  $H_s^p(w)$ consists of continuous functions on $\overline{\B}^n$, and consequently, the Carleson measures are just the finite ones.  If $\tau-sp<1$, we show that the Carleson measures for $H_s^p(w)$ and $K_s[L^p(w)]$ coincide, whereas  if $\tau-sp\geq 1$, this coincidence may fail.

   As it happens in the unweighted case (see  \cite{cohnverbitsky2}), the proof of the characterization of the Carleson measures for $H_s^p(w)$ will be based in the construction of   weighted holomorphic potentials,  with control of their $H_s^p(w)$-norm. In fact, technical reasons give that it is  convenient to deal with weighted Triebel-Lizorkin spaces which, on the other hand,   have   interest on their own. In the second section we study   these spaces.  
    If $s\geq 0$, we will write $[s]^+$ the integer part of $s$ plus $1$. Let $1<p<+\infty$, $1\leq q\leq +\infty$, and $s\geq 0$. The weighted holomorphic Triebel-Lizorkin space $HF_s^{pq}(w)$ when $q<+\infty$ is the space of holomorphic functions $f$ in $\B^n$ for which
\begin{equation*}\begin{split}&||f||_{HF_s^{pq}(w)}=\\&\left( \int_{{\bf S}^n}\left( \int_0^1 |((I+R)^{[s]^+}f)(r\zeta)|^q(1-r^2)^{([s]^+-s)q-1}dr\right)^{\frac{p}{q}}w(\zeta) d\sigma(\zeta)\right)^\frac{1}{p}<+\infty,
\end{split}\end{equation*}  
 whereas when $q=+\infty$, 
  $$||f||_{HF_s^{p\infty}(w)}=\left(\int_{{\bf S}^n}\left(\sup_{0<r<1}|(( I+R)^{[s]^+}f)(r\zeta)|(1-r^2)^{[s]^+-s}\right)^pw(\zeta) d\sigma(\zeta)\right)^\frac1{p}<+\infty,
$$
where $I$ denotes the identity operator.

The  Section $2$ is devoted  to the general theory of weighted holomorphic Triebel-Lizorkin spaces.  We give   different equivalent definitions of the spaces $HF_s^{p,q}(w)$ in terms of admissible area functions, we give duality theorems on these spaces, we study some relations of inclusion among them 
and we also obtain that when $q=2$, the weighted Triebel-Lizorkin space $HF_s^{p2}(w)$ coincides with the weighted Hardy-Sobolev space $H_s^p(w)$.

The main result in Section 3 is the characterization of the Carleson measures for $H_s^p(w)$, when $0<\tau-sp<1$, in terms of a positive kernel.
\begin{theoremC}
Let $1< p<+\infty$, $w$  an $A_p$-weight, and $\mu$ a finite positive Borel measure on $\B^n$. Assume that $w$ is doubling of order $\tau$, for some $\tau<1+sp$. We then have that the following statements are equivalent:
 
 (i) $||K_s(f)||_{L^p(d\mu)}\leq C||f||_{ L^p(w)}$.
 
 (ii) $||f||_{L^p(d\mu)}\leq C|| f||_{H_s^p(w)}$.
\end{theoremC}
The proof relies on the construction of weighted holomorphic potentials, with control of their weighted Hardy-Sobolev norm.

We also gives examples of the sharpness of the above theorem. We show that if $p=2$ and $\tau>1+sp$, $n<\tau<n+1$, then there exists $w$ in $A_2\cap D_\tau$ and a measure $\mu$ on ${\bf S}^n$ which is Carleson for $H_s^2(w)$, but it is not Carleson for $K_s[L^2(w)]$.

 Finally, the usual remark on  notation: we will adopt the  convention
of using the same letter for various absolute constants whose
values may change in each occurrence, and we will write $A\preceq
B$ if there exists an absolute constant $M$ such that $A\leq MB$.
 We will say that two quantities $A$ and $B$ are equivalent if both
$A\preceq B$ and $B\preceq A$, and, in that case, we will write
$A\simeq B$.

\section{Weighted holomorphic Triebel-Lizorkin spaces }
\label{section1}
In this section we will introduce weighted holomorphic Triebel-Lizorkin spaces, and we will obtain characterizations in terms of Littlewood-Paley functions and admissible area functions. These characterizations, known in the unweighted case, will be used in the following sections. 

We begin recalling  some simple facts about $A_p$ weights that we will need later.
It is well known that $A_\infty=\bigcup_{1<p<+\infty}A_p$ and that any $A_p$ weight satisfies a doubling condition. We recall that a weight $w$ satisfies a doubling condition of order $\tau$, $\tau>0$, if there exists $C>0$,   such that for any nonisotropic ball $B\subset {\bf S}^n$, and any $k\geq 0$, $W(2^kB)\leq C2^{\tau k}W(B)$. We will say that this weight    $w$ is in $D_\tau$.  In fact, if $w\in A_p$, there exists $p_1<p$ such that $w$ is also in $A_{p_1}$, and consequently we have that $w\in D_{\tau}$ for $\tau=np_1<np$, (see \cite{strombergtorchinsky}).  

Examples of $A_p$ weights    can be obtained as follows: if $\zeta=(\zeta',\zeta_n)$, and $w(\zeta)=(1-|\zeta'|^2)^\varepsilon$, we then have that $w\in A_p$ if $-1<\varepsilon<p-1$. We also have that for this weight, $w\in D_\tau$, $\tau=n+\varepsilon$.

The following lemma gives the natural relationships between the spaces  $L^p(w)$, $w\in A_p$, and the Lebesgue spaces $L^q(d\sigma)$.
\begin{lemma}\label{weightedspaces}
Let $1<p<+\infty$, and $w$ be an $A_p$-weight. We then have:

(i) There exists $1<p_1<p$ such that $L^p(w)\subset L^{p_1}(d\sigma)$.

(ii) There exists $ p_2>p$ such that $L^{p_2}(d\sigma) \subset L^{p }(w)$.
\end{lemma}

We now proceed to study the weighted holomorphic Triebel-Lizorkin spaces $H_s^{pq}(w)$ already defined in the introduction. We begin with some definitions.
If $1<q\leq +\infty$, $k$ an integer such that $k>s\geq0$,  and $\zeta\in{\bf S}^n$, the Littlewood-Paley type functions are given by
$$
A_{1,k,q,s}(f)(\zeta)=\left(\int_0^1 |(I+R)^kf(r\zeta)|^q(1-r^2)^{(k-s)q-1}dr \right)^\frac1{q},
$$
when $q<+\infty$, and
$$
A_{1,k,\infty,s}(f)(\zeta)=\sup_{0<r<1}|(I+R)^kf(r\zeta)|(1-r^2)^{k-s},$$
when $q=+\infty$.

If $\alpha>1$, $\zeta\in {\bf S}^n$, we denote by $D_\alpha(\zeta)$, $\alpha>1$   the admissible region given by $D_\alpha(\zeta)=\{z\in\B^n\,;\, |1-z\overline{\zeta}|<\alpha(1-|z|)\,\}$. We   introduce the admissible area function
$$
A_{\alpha,k,q,s}(f)(\zeta)=\left(\int_{D_\alpha(\zeta)} |(I+R)^kf(z)|^q( 1-|z|^2)^{(k-s)q-n-1}dv(z)\right)^\frac1{q},
$$
when $q<+\infty$, where $dv$ is the Lebesgue measure on $\B^n$, and in case $q=+\infty$, 
$$
A_{\alpha,k,\infty,s}(f)(\zeta)=\sup_{z\in D_\alpha(\zeta)}|(I+R)^kf(z)|(1-|z|^2)^{k-s},$$
 when $q=+\infty$.

 Our first goal is to obtain that if $1<p<+\infty$, $1<q<+\infty$ and  $w$ is an $A_p$ weight, then an holomorphic function is in $HF_{s}^{p,q}(w)$ if and only if $A_{\alpha,k,q,s}(f)\in L^p(w)$, for some (and then for all) $\alpha\geq 1$ and
$k>s$.
We will follow  the ideas in \cite{ortegafabrega}. For the sake of completeness, we will sketch  the modifications needed to obtain the weighted case.

If $1<p<+\infty$, $1<q\leq+\infty$ we denote by $\displaystyle{ L^p(w)(L^q_1)= L^p(w)(L^q (\frac{2nr^{2n-1}}{1-r^2}}dr))$ the mixed-norm space of measurable functions $f$ in ${\bf S}^n \times [0,1]$ such that
$$||f||_{p,q,w}=\left( \int_{{\bf S}^n} \left( \int_0^1 |f(r\zeta)|^q \frac{2nr^{2n-1}}{1-r^2}dr\right)^\frac{p}{q} w(\zeta) d\sigma(\zeta)\right)^\frac1{p}<+\infty.$$
Also if $\alpha>1$, and $E_\alpha(z)=\left( \int_{{\bf S}^n} \chi_{D_\alpha(\zeta)}(z)d\sigma(\zeta)\right)^{-1}\simeq (1-|z|^2)^{-n}$, we denote by $L^p(w)(L^q_\alpha)$   the mixed-norm space of measurable functions $f$ defined in  
${\bf S}^n \times \B^n $
such that
$$||f||_{\alpha,p,q,w}=\left(\int_{{\bf S}^n} \left( \int_{\B^n} |f(z,\zeta)|^q 
\frac{E_\alpha(z)}{(1-|z|^2)}dv(z)\right)^\frac{p}{q} w(\zeta)d\sigma(\zeta)\right)^\frac1{p}<+\infty.$$

We denote by $F^{\alpha,\,p,q}(w)$ the space of measurable functions on $\B^n$ such that $$J_\alpha f(\zeta,z)= \chi_{D_\alpha(\zeta)}(z)f(z)$$ is in $ L^p(w)(L^q_\alpha)$, normed with the norm induced by $ ||\cdot||_{\alpha,p,q,w}$. We also introduce the space $F^{1,p,q}(w)$ of measurable functions on $\B^n$ such that $J_1f(\zeta,r)= f(r\zeta)$ is in $ L^p(w)(L^q_1)$. 

The representation of the dual of a mixed-norm space, see \cite{berglofstrom}, gives that if $1<p,q<+\infty$, the dual space of $ L^p(w)(L^q_1)$ is $L^{p'}(w)(L^{q'}_1)$, $1/p+1/p'=1$, $1/q+1/q'=1$, and that if $f\in F^{1,p,q}(w)$, $g\in F^{1,p',q'}(w)$ the pairing is given by
\begin{equation*} 
(f,g)=\int_{{\bf S}^n}\left(\int_0^1 f(r\zeta) \overline{g(r\zeta)}\frac{2nr^{2n-1}}{ 1-r^2} dr \right)w(\zeta) d\sigma(\zeta).
 \end{equation*}
Analogously, the dual space of $ L^p(w)(L^q_\alpha)$ is $L^{p'}(w)(L^{q'}_\alpha)$, and   if $f\in F^{\alpha,\, p,q}(w)$, $g\in F^{\alpha,\, p',q'}(w)$ the pairing is given by
\begin{equation*}\begin{split} (f,g)_\alpha= &\int_{\B^n}\int_{{\bf S}^n}f(z)\overline{g(z)} \chi_{D_\alpha(\zeta)}(z) w(\zeta) d\sigma(\zeta) \frac{dv(z)}{(1-|z|^2)^{n+1}}  \\=& \int_{\B^n}f(z)\overline{g(z)}\frac{E_\alpha^w(z)}{(1-|z|^2)^{n+1}}dv(z),
\end{split}\end{equation*}
where  $E_\alpha^w(z)=\int_{{\bf S}^n} \chi_{D_\alpha(\zeta)}(z) w(\zeta) d\sigma(\zeta) $. 

Observe that if   we write $z_0=\frac{z}{|z|}$, the doubling property of $w$ gives that $E_\alpha^w(z)\simeq W(B(z_0,(1-|z|)))$. From now on we will write $B_z=B(z_0,(1-|z|))$.

 We begin with two lemmas that are weighted versions of 
 Lemmas 2.2. and 2.3 in \cite{ortegafabrega}, and whose proofs we omit. We recall that if $\psi$ is a measurable function on ${\bf S}^n$, the weighted Hardy-Littlewood maximal function is given by
 $$
 M_{HL}^w(\psi)( \zeta)=\sup_{B\ni \zeta}\frac1{W(B)}\int_B |\psi(\eta)|w(\eta)d\sigma (\eta).$$

 \begin{lemma} \label{lemma2.1}
There exist $C>0$, $N_0>0$ such that for any $z\in D_\alpha(\zeta)$, $N\geq N_0$,
$$
\frac{(1- |z|^2)^{n+N}}{W(B_z)}\int_{{\bf S}^n} \frac{|\psi(\eta)|}{|1-z\overline{\eta}|^{n+N}}w(\eta)d\sigma(\eta)\leq C M_{HL}^w(\psi)( \zeta).
$$
 
 \end{lemma}

 \begin{lemma} \label{lemma2.2}
Let $\alpha>1$. There exists $C>0$, such that for any $z\in D_\alpha(\zeta)$,  
$$
\frac{1}{W(B_z)}\int_{{\bf S}^n}  \chi_{D_\alpha(\eta)}(z)|\psi(\eta)| w(\eta)d\sigma(\eta)\leq C M_{HL}^w(\psi)( \zeta).
$$
 
 \end{lemma}

 \begin{theorem}\label{theorem2.3}
 Let $1<p<+\infty$, $1\leq q\leq+\infty$, and $\alpha\geq 1$. Then the space $F^{\alpha,\, p,q}(w)$ is a retract of $ L^p(w)(L_\alpha^q)$.
 \end{theorem}
{\bf Proof of Theorem \ref{theorem2.3}:}\par
The fact that $J_1$ is an isometry between $F^{1,\, p,q}(w)$ and $ L^p(w)(L_1^q)$ gives the theorem for the case $\alpha=1$.

If $\alpha>1$, we introduce the averaging operator
$$
A_\alpha(\varphi)(z)=\frac{1}{E_\alpha^w(z)}\int_{{\bf S}^n} \chi_{D_\alpha(\eta)}(z)\varphi(\eta, z)w(\eta) d\sigma(\eta).
$$
The definition of $E_\alpha^w(z)$ gives that $A_\alpha\circ J_\alpha$ is the identity operator on $F^{\alpha,\, p,q}(w)$. So, in order to finish the theorem, we need to show that $A_\alpha$ maps $ L^p(w)(L^q_\alpha)$ to $F^{\alpha,\,p,q}(w)$.
We consider first the case $1\leq q\leq p<+\infty$. Let $m=\frac{p}{q}\geq 1$ and let $m'$ be the conjugate exponent of $m$. We then have by duality that
$$||A_\alpha(\varphi)||_{\alpha,p,q,w}^q=\sup_{ ||\psi||_{L^{m'}(w)}\leq 1} |\int_{{\bf S}^n} \int_{D_\alpha(\zeta)}|A_\alpha(\varphi)(z)|^q\frac{dv(z)}{(1-|z|^2)^{n+1}} \psi(\zeta)w(\zeta) d\sigma(\zeta)|.$$
 
Now H\"older's inequality gives that
$$|A_\alpha(\varphi)(z)|^q\leq\frac1{E_\alpha^w(z)}\int_{{\bf S}^n} |\varphi(\eta,z)|^q \chi_{D_\alpha(\eta)}(z)w(\eta)d\sigma(\eta).$$ 
Hence, by Lemma \ref{lemma2.2}
\begin{equation*}\begin{split}
&||A_\alpha(\varphi)||_{\alpha,p,q,w}^q  \\&\preceq
 \sup_{ ||\psi||_{L^{m'}(w)}\leq 1}\int_{{\bf S}^n} \int_{\B^n}\frac1{E_\alpha^w(z)}\chi_{D_\alpha(\zeta)}(z) \int_{{\bf S}^n} \chi_{D_\alpha(\eta)}(z)|\varphi(\eta,z)|^q w(\eta) d\sigma(\eta)\\
 &\times \frac{dv(z)}{(1-|z|^2)^{n+1}}|\psi(\zeta)|w(\zeta)d\sigma(\zeta)\\
 &\preceq\sup_{||\psi||_{L^{m'}(w)}\leq 1}\int_{{\bf S}^n} \int_{\B^n} 
 |\varphi(\eta,z)|^q\frac{dv(z)}{(1-|z|^2)^{n+1}} w(\eta)M_{HL}^w(\psi)(\eta)d\sigma(\eta),
\end{split}\end{equation*}
 Next, H\"older's inequality with exponent $m=\frac{p}{q}$ gives that the above is bounded by
\begin{equation*}\begin{split}
&\sup_{ ||\psi||_{L^{m'}(w)}\leq 1}||M_{HL}^w\psi||_{L^{m'}(w)}\left( \int_{{\bf S}^n}\left( \int_{\B^n}   |\varphi(\eta,z)|^q\frac{dv(z)}{(1-|z|^2)^{n+1}} \right)^\frac{p}{q} w(\eta) d\sigma(\eta) \right)^\frac{q}{p}\\
&\leq\sup_{ ||\psi||_{L^{m'}(w)}\leq 1}||\psi||_{L^{m'}(w)} ||\varphi||_{\alpha,\,p,q,w}^q,
\end{split}\end{equation*}
where we have used that since $w$ is a doubling measure, the weighted Hardy-Littlewood maximal function is bounded from $L^{m'}(w)$ to ${L^{m'}(w)}$.
That finishes the proof of the theorem when $q\leq p$.

So we are lead to deal with the case $1<p< q\leq +\infty$, which can be easily obtained from the previous case using the duality in the mixed-norm spaces $ L^p(w)(L^q_\alpha)$.\qed

This result can be used as in the unweighted case to obtain a characterization of the dual spaces of the  weighted spaces $F^{\alpha,\,p,q}(w)$.
\begin{corollary}\label{corollary2.4}
Let $1<p<+\infty$, $1< q<+\infty$, $\alpha>1$, and $w$ an $A_p$-weight. Then the dual of $F^{\alpha,\,p,q}(w)$ is $F^{\alpha,\,p',q'}(w)$ with the pairing given by $(f,g)_\alpha$.
\end{corollary}

The following proposition will be needed in the proof of the main theorem in this section. If $N>0$, $M>0$, we consider the operators  defined by
$$
P^{N,M}f(y)=\int_{\B^n} f(z)\frac{(1-|z|^2)^N (1-|y|^2)^M}{|1-z\overline{y}|^{n+1+N+M}}dv(z),\,\,\,  y\in\B^n .
$$

\begin{theorem}\label{theorem2.5.1}
Let $1<p<+\infty$, $1\leq q< +\infty$, $\alpha,\beta\geq 1$, and $w$ an  $A_p$ weight. Then there exists $N_0>0$ such that for any $N\geq N_0$ and any $M>0$, the operator $P^{N,M}$ is continuous from $F^{\alpha,\,p,q}(w)$ to $F^{\beta,\,p,q}(w)$.
\end{theorem}
 
{\bf Proof of Theorem \ref{theorem2.5.1}:}\par

We begin with the case $\alpha,\beta>1$. The case where $1\leq q\leq p<+\infty$ can be deduced following the scheme of \cite{ortegafabrega}, using Lemma \ref{lemma2.1}.
 
In the case $1<p<q< +\infty$ we apply duality in the mixed norm space and obtain
\begin{equation}\begin{split}\label{equation2.2}
||P^{N,M}(f)||_{\beta,\,p,q,w}^q&=\sup_{  ||g||_{\beta,p'q',w}\leq 1}|\int_{\B^n} P^{N,M}(f)(y)\overline{g(y)} \frac{E_\beta^w(y)}{(1-|y|^2)^{n+1}}dv(y)|\\
& \leq\sup_{  ||g||_{\beta,p'q',w}\leq 1} (f, \widetilde{P}^{M-1,N+1}(g))_\alpha,
\end{split}\end{equation}
  where
\begin{equation}\label{equation2.3}
\widetilde{P}^{R,S}(g)(z)=  \int_{\B^n} \frac{(1-|y|^2)^R(1-|z|^2)^Sg(y)}{|1-y\overline{z}|^{n+1+R+S}} \frac{ E_\beta^w(y)}{(1-|y|^2)^n}\frac{(1-|z|^2)^n}{E_\alpha^w(z)}dv(y).
\end{equation}
Observe that when $w\equiv 1$, then $\widetilde{P}^{M,N}(f)\simeq P^{M,N}(f)$.
Here we are led to obtain that the operator $\widetilde{P}^{M-1,N+1}$ maps boundedly $F^{\beta,p',q'}$ to $F^{\alpha,p',q'}$,  provided $p< q$. 
If we claim this proposition, we finish the proof of the theorem.
Using (\ref{equation2.2}), and applying H\"older's inequality,
\begin{equation*}\begin{split}
& ||P^{N,M}(f)||_{\beta,\,p,q,w}^q= \sup_{ ||g||_{\alpha,p'q',w}\leq 1 }(f, \widetilde{P}^{M-1,N-1}(g))_\alpha\leq\\
&\sup_{ ||g||_{\alpha,p'q',w}\leq 1 }||f||_{\alpha,p,q,w}||\widetilde{P}^{M-1,N-1}(g)||_{\alpha, p',q',w}\leq C\sup||f||_{\alpha,p,q,w} .
\end{split}\end{equation*}
  The cases $\alpha=1$ and $\beta=1$ are proved in a simmilar way. 
  
  To finish the theorem we will prove the claim. Changing the notation, it is enough to prove:
\begin{proposition}\label{proposition2.6}
Let $1<q<p<+\infty$,   $\alpha,\beta\geq 1$, and $w$ an $A_p$ weight. We then have that there exists $N_0>0$ such that for any $N\geq N_0$ and any $M\geq 0$,

(i) $\widetilde{P}^{M,N}(1)<+\infty$.

(ii)  The operator $\widetilde{P}^{M,N}$ is continuous from $F^{\alpha,\,p,q}(w)$ to $F^{\beta,\,p,q}(w)$.
\end{proposition}

{\bf Proof of Proposition \ref{proposition2.6}:}\par
Let us begin with (i). From the definition of $E_\alpha^w(z)$ and Fubini's theorem,
\begin{equation*}\begin{split}
& \int_{\B^n} \frac{(1-|z|^2)^M}{|1-z\overline{y}|^{n+1+M+N}}\frac{E_\alpha^w(z)}{(1-|z|^2)^n}dv(z)\\
&=\int_{{\bf S}^n} \int_{D_\alpha(z)}\frac{(1-|z|^2)^M}{|1-z\overline{y}|^{n+1+M+N}}\frac{dv(z)}{ (1-|z|^2)^n}w(\zeta) d\sigma(\zeta) \preceq 
 \int_{{\bf S}^n} \frac1{|1-y\overline{\zeta}|^{n+N}}w(\zeta)d\sigma(\zeta),
\end{split}\end{equation*}
where in last inequality we have used Lemma 2.7 in \cite{ortegafabrega} since $M>-1$.
 
 Next, let $B_k=B(y_0, 2^k(1-|y|^2))$, $k\geq 0$, where $y_0=\frac{y}{|y|}$. Since $w$ is doubling and $E_\alpha^w(y)\simeq W(B_0)$ give  that $W(B_k)\leq C^k E_\alpha^w(y)$. Consequently
 \begin{equation*}\begin{split}
 &\int_{{\bf S}^n} \frac1{|1-y\overline{\zeta}|^{n+N}}w(\zeta)d\sigma(\zeta)\\
 &\preceq\sum_k \int_{B_k}\frac{w(\zeta)d\sigma(\zeta)}{(2^k(1-|y|^2))^{n+N}}\preceq \frac{E_\alpha^w(y)}{(1-|y|^2)^{n+N}} \sum_k \frac{C^k}{2^{k(n+N)}}\preceq\frac{E_\alpha^w(y)}{(1-|y|^2)^{n+N}},
 \end{split}\end{equation*}
if $N$ is chosen sufficiently large. That finishes the proof of (i).

Since $m=\frac{p}{q}> 1$, duality gives that
\begin{equation}\label{equation2.1}
||\widetilde{P}^{M,N}(f)||_{\beta,\,p,q,w}^q=\sup_{   ||\psi||_{L^{m'}(w)}\leq 1}|\int_{{\bf S}^n} \int_{D_\beta(\zeta)}|\widetilde{P}^{M,N}f(y)|^q \frac{dv(y)}{(1-|y|^2)^{n+1}}\overline{\psi(\zeta)}w(\zeta) d\sigma(\zeta)|.
\end{equation}
 
 Next, H\"older's inequality shows that if $0<\varepsilon<N$ then
\begin{equation*}\begin{split}
&|\widetilde{P}^{M,N}(f)(y)|^q\leq \int_{\B^n} |f(z)|^q\frac{(1-|z|^2)^M(1-|y|^2)^{N-\varepsilon}}{|1-z\overline{y}|^{n+1+M+N-\varepsilon}}\frac{E_\alpha^w(z)}{(1-|z|^2)^n}\frac{ (1-|y|^2)^n}{E_\alpha^w(y)}dv(z)\\&\times \left( \int_{\B^n} \frac{(1-|z|^2)^M(1-|y|^2)^{N+\varepsilon \frac{q'}{q}}}{|1-z\overline{y}|^{n+1+M+N+\varepsilon \frac{q'}{q}}}\frac{E_\alpha^w(z)}{(1-|z|^2)^n}\frac{ (1-|y|^2)^n}{E_\alpha^w(y)}dv(z) \right)^\frac{q}{q'}\\&\preceq
\int_{\B^n} |f(z)|^q\frac{(1-|z|^2)^M(1-|y|^2)^{N-\varepsilon}}{|1-z\overline{y}|^{n+1+N+M-\varepsilon}}\frac{E_\alpha^w(z)}{(1-|z|^2)^n}\frac{ (1-|y|^2)^n}{E_\alpha^w(y)}dv(z),
\end{split}\end{equation*}
where in last inequality we have used (i).

Consequently,
\begin{equation}\begin{split} \label{equation2.4}
&||\widetilde{P}^{M,N}(f)||_{\beta,\,p,q,w}^q\leq C\sup_{   ||\psi||_{L^{m'}(w)}\leq 1}|\int_{{\bf S}^n} \int_{y\in D_\beta(\zeta)}\int_{\B^n} \frac{|f(z)|^q (1-|z|^2)^M(1-|y|^2)^{N-\varepsilon}}{|1-z\overline{y}|^{n+1+N+M-\varepsilon}}\\
&\times \frac{E_\alpha^w(z)}{(1-|z|^2)^{n }}\frac{(1-|y|^2)^{n }}{E_\alpha^w(y)}
dv(z)\frac{dv(y)}{(1-|y|^2)^{n+1}}\psi(\zeta)w(\zeta)d\sigma(\zeta)|\\
&=C \sup_{  ||\psi||_{L^{m'}(w)}\leq 1}|\int_{{\bf S}^n}\int_{\B^n} \int_{D_\beta(\zeta)}
\frac{(1-|y|^2)^{N+n-\varepsilon}}{|1-z\overline{y}|^{n+1+N+M-\varepsilon}} \frac{dv(y)}{E_\alpha^w(y)(1-|y|^2)^{n+1}}\\
&\times  |f(z)|^q (1-|z|^2)^{M-n} E_\alpha^w(z) dv(z) \psi(\zeta)w(\zeta) d\sigma(\zeta)|.
\end{split}\end{equation}

Next, if $y\in D_\beta(\zeta)$, $E_\alpha^w(y)\simeq W(B_y)\simeq W(B(\zeta, (1-|y|^2))$, and $|1-z\overline{y}|\simeq (1-|y|^2)+|1-z\overline{\zeta}|$.  

Assume first that $|1-z\overline{\zeta}|\leq 1$. Hence, 
\begin{equation}\begin{split} \label{equation2.5} 
&\int_{D_\beta(\zeta)}
\frac{(1-|y|^2)^{N+n-\varepsilon}}{|1-z\overline{y}|^{n+1+N+M-\varepsilon}} \frac{dv(y)}{E_\alpha^w(y)(1-|y|^2)^{n+1}} \\
&\simeq\int_{\B^n} \frac{(1-|y|^2)^{N-\varepsilon}}{((1-|y|^2)+|1-z \overline{\zeta}|)^{n+1+N+M-\varepsilon}}\chi_{D_\beta(\zeta)}(y) \frac{(1-|y|^2)^n}{W(B(\zeta, 1-|y|^2))}\frac{dv(y)}{(1-|y|^2)^{n+1}},\end{split}\end{equation}
which by integration in polar coordinates 
\begin{equation*}\begin{split}&\int_0^1\frac{r^{2n-1} (1-r^2)^{N+n-\varepsilon}}{((1-r^2)+|1-z\overline{\zeta}|)^{n+1+N+M-\varepsilon}}  \frac{dr}{(1-r^2)W(B(\zeta,C(1-r^2)))}  \\& \simeq
  \int_0^{|1-z\overline{\zeta}|} \frac{t^{N+n-\varepsilon-1}}{(t+|1-z\overline{\zeta}|)^{n+1+N+ M-\varepsilon} }\frac{dt}{W(B(\zeta,t))}\\&+\int_{|1-z\overline{\zeta}|}^1 \frac{t^{N+n-\varepsilon-1}}{(t+|1-z\overline{\zeta}|)^{n+1+N+ M-\varepsilon}}\frac{dt}{W(B(\zeta,t))}=I+II.
\end{split}\end{equation*}
In $I$ we have that $(t+|1-z\overline{\zeta}|)\simeq |1-z\overline{\zeta}|$, and, since $w\in A_p$, $$\frac{t^n}{W(B(\zeta,t))}\simeq \left( \frac1{t^n} \int_{B(\zeta,t)}w^{-(p'-1)}\right)^{p-1}.$$ Thus we obtain  that
\begin{equation*}\begin{split}
&I\simeq
\int_0^{|1-z\overline{\zeta}|} \frac{t^{N-\varepsilon-1}}{ |1-z\overline{\zeta}|^{n+1+N+ M-\varepsilon}}\left( \frac1{t^n} \int_{B(\zeta,t)}w^{-(p'-1)}\right)^{p-1}dt
\\&\preceq \left(\int_{B(\zeta,|1-z\overline{\zeta}|)}w^{-(p'-1)}\right)^{p-1} \frac{1}{|1-z\overline{\zeta}|^{n+1+N+M-\varepsilon}}
\int_0^{|1-z\overline{\zeta}|}t^{N-\varepsilon-n(p'-1)-1} dt\\ 
&\preceq \frac{1}{|1-z\overline{\zeta}|^{M+1}}\frac1{W(B(z_0, |1-z\overline{\zeta}|))},\end{split}\end{equation*}
where  we have used that $N>0$ is chosen big enough, and that $w$ satisfies the $A_p$ condition.

In II, $(t+|1-z\overline{\zeta}|)\simeq t$, and since $M+1>0$, we have
\begin{equation*}\begin{split}&II\simeq \int_{|1-z\overline{\zeta}|}^1 \frac{1}{t^{M+2}}\frac{dt}{W(B(\zeta,t))}\\&\leq
\int_{|1-z\overline{\zeta}|}^1 \frac{1}{t^{M+2}}\frac{dt}{W(B(\zeta,|1-z\overline{\zeta}|))}\preceq  \frac{1}{|1-z\overline{\zeta}|^{M+1}}\frac{1}{W( B(z_0,|1-z\overline{\zeta}|))}.
\end{split}\end{equation*}
If $|1-z\overline{\zeta}|>1$, then we have that $(1-r^2)+|1-z\overline{\zeta}|\simeq 1$. We return to  (\ref{equation2.5}) and obtain
\begin{equation*}\begin{split}&\int_0^1 \frac{(1-r^2)^{N+n-\varepsilon-1}dr}{((1-r)^2+|1-z\overline{\zeta} |)^{n+1+N+M-\varepsilon}W(B(\zeta,1-r^2))} \\
& \preceq \left(\int_{B(\zeta,1)}w^{-\frac{p'}{p}}\right)^{\frac{p}{p'}}\int_0^1 t^{N-\varepsilon-n\frac{p}{p'}-1}dt\preceq
\frac{1}{|1-z\overline{\zeta}|^{M+1}}\frac{1}{W( B(z_0,|1-z\overline{\zeta}|))}.
\end{split}\end{equation*}

Then we have in any case that (\ref{equation2.5}) is bounded by
  $$\frac{1}{|1-z\overline{\zeta}|^{M+1}}\frac{1}{W( B(z_0,|1-z\overline{\zeta}|))}.$$ In consequence, we return to (\ref{equation2.4}) and we obtain  
\begin{equation}\begin{split}\label{equation2.6}
&||\widetilde{P}^{M,N}(f)||_{\beta,\,p,q,w}^q \\&\preceq  
\sup_{ ||\psi||_{L^{m'}(w)}\leq 1} |\int_{{\bf S}^n}\int_{\B^n}   |f(z)|^q  \frac{(1-|z|^2)^{M-n}   E_\alpha^w(z)}{|1-z\overline{\zeta}|^{M+1} W(B(z_0,|1-z\overline{\zeta}|))} \psi(\zeta)dv(z)w(\zeta) d\sigma(\zeta)|\\&\preceq  \sup_{ ||\psi||_{L^{m'}(w)}\leq 1} |\int_{{\bf S}^n}\int_{\B^n}   |f(z)|^q   (1-|z|^2)^{M-n} \chi_{D_\alpha(\eta)}(z) \int_{{\bf S}^n} \frac{\psi(\zeta)w(\zeta) d\sigma(\zeta)}{|1-z\overline{\zeta}|^{M+1} W(B(z_0,|1-z\overline{\zeta}|))}\\&dv(z) w(\eta) d\sigma(\eta)|.
\end{split}\end{equation}
Next, if $z\in D_\alpha(\eta)$, $B(\eta, |1-z\overline{\zeta}|)\subset B(z_0, C|1-z\overline{\zeta}|)$, and if $B_k=B(\eta, 2^k(1-|z|^2))$, $k\geq 0$    and $\zeta\in B_{k+1}\setminus B_k$, $|1-z\overline{\zeta}|\simeq 2^k(1-|z|^2)$. Thus
\begin{equation*}\begin{split}
&\int_{{\bf S}^n} \frac{|\psi(\zeta)|w(\zeta) d\sigma(\zeta)}{|1-z\overline{\zeta}|^{M+1} W(B(z_0,|1-z\overline{\zeta}|))} \\&\preceq
 \frac1{(1-|z|^2)^{M+1} W(B(\eta, 1-|z|^2))}\int_{B_0} |\psi(\zeta)|w(\zeta) d\sigma(\zeta) \\&+
\sum_{k\geq 1} \frac1{2^{k(M+1)}(1-|z|^2)^{M+1} W(B(\eta, 2^k(1-|z|^2)))}\int_{B_k} |\psi(\zeta)|w(\zeta) d\sigma(\zeta)\\&\preceq
\frac1{(1-|z|^2)^{M+1}}\sum_{k\geq 0}\frac1{2^{k(M+1)}} M_{HL}^w(\psi)(\eta)\preceq \frac1{(1-|z|^2)^{M+1}}  M_{HL}^w(\psi)(\eta).
\end{split}\end{equation*}
Plugging the above estimate in (\ref{equation2.6}) and using H\"older's inequality with exponent $m=\frac{p}{q}$, we obtain
\begin{equation*}\begin{split} 
&||\widetilde{P}^{M,N}(f)||_{\beta,\,p,q,w}^q \\&  \preceq\sup_{\psi\in L^{m'}(w)}  \int_{{\bf S}^n}\int_{\B^n}   |f(z)|^q  \frac1{(1-|z|^2)^{n+1}}\chi_{D_\alpha(\eta)}(z) dv(z)M_{HL}^w(\psi)(\eta)w(\eta) d\sigma(\eta) \\
&\preceq \sup_{\psi\in L^{m'}(w)}  ||f||_{\alpha,p,q,w}^q||M_{HL}^w(\psi)||_{L^{m'}(w)}^q \preceq ||f||_{\alpha,p,q,w}^q.\qed
\end{split}\end{equation*}

 We deduce from the previous theorem  the following characterization of the weighted holomorphic Triebel-Lizorkin spaces.  If $f\in H(\B^n)$, $s,t>0$, let $$L_s^tf(z)=(1-|z|^2)^{t-s} (I+R)^tf(z).$$
\begin{theorem}\label{theorem2.5}
Let $1<p<+\infty$, $1<q<+\infty$,  $t>s\geq 0$ and $\alpha\geq 1$. Let 
$$HF_s^{\alpha,\,t,\,p,q}(w)=\{ f\in H(\B^n);\, ||L_s^tf||_{\alpha,p,q}<+\infty\}.$$
 Then $HF_s^{\alpha,\,t,\,p,q}(w)=HF_s^{pq}(w)$.
\end{theorem}
{\bf Proof of theorem \ref{theorem2.5}:}\par
If $s<t_0<t_1$, $\alpha,\beta\geq 1$, we just need to check that $HF_s^{\alpha,\,t_0,\,p,q}(w)=HF_s^{\beta,\,t_1,\,p,q}(w)$. Any holomorphic function $f$ on $\B^n$ satisfying that  $L_s^tf(z)\in F^{\alpha,\,p,q}(w)$ is in $A^{-\infty}(\B^n)$, the space of holomorphic functions in $\B^n$ for which there exists $k>0$ such that $\sup_z(1-|z|^2)^k|f(z)|<+\infty$. Consequently, $f$ and its derivatives have a representation formula via the reproducing kernel $c_N\frac{(1-|z|^2)^N}{(1-\overline{z}y)^{n+1+N}}$, for $N>0$ sufficiently large and an adequate constant $c_N$. Once we have made this observation, we  can reproduce the arguments in \cite{ortegafabrega} and obtain
$$(I+R)^{t_0}f(y)=C_N\int_{\B^n}(I+R)^{t_1}f(z)(I+R_y)^{t_0-t_1} \frac{(1-|z|^2)^N}{(1-y\overline{z})^{n+1+N}}dv(z).$$
Since for $m>0$ we have that
\begin{equation}\label{novalabel1}(I+R)^{-m}g(y)=\frac1{\Gamma(m)}\int_0^1 \left( \log \frac1{r}\right)^{m-1}g(ry)dr,\end{equation}
we obtain
\begin{equation*}\begin{split}
||L_s^{t_0}f||_{\alpha,p,q,w}&\preceq || \int_{\B^n} |(I+R)^{t_1}f(z)|\frac{(1-|z|^2)^N(1-|y|^2)^{t_0-s}}{|1-\overline{z} y|^{n+1+N+t_0-t_1}}dv(z)||_{\alpha,p,q,w}\\&=||P^{N-t_1+s,t_0-s}(|L_s^{t_1}f|)||_{\alpha,p,q,w},\end{split}\end{equation*}
and we just have to apply Theorem \ref{theorem2.5.1} to finish the proof.\qed 
\begin{theorem}\label{theorem2.7}
 Let $1<p<+\infty$, $1<q<+\infty$,  $w$   an $A_p$-weight, and $f$ a holomorphic function. Then the following assertions are equivalent:
 
 (i) $f$ is in $HF_s^{pq}(w)$.
 
 (ii)   $A_{\alpha,k,q,s}(f)\in L^p(w)$, for some   $\alpha\geq 1$ and
$k>s$.

(iii) $A_{\alpha,k,q,s}(f)\in L^p(w)$, for all   $\alpha\geq 1$ and
$k>s$.
\end{theorem}
 
Our next result studies some inclusion relationships between different weighted holomorphic Triebel-Lizorkin spaces.
\begin{theorem}\label{theorem2.8}
Let $1<p<+\infty$, $1\leq q_0\leq q_1\leq +\infty$,  $s\geq 0$ and let $w$ be an $A_p$-weight. We then have 
$$HF_s^{pq_0}(w)\subset HF_s^{pq_1}(w).$$
\end{theorem}
{\bf Proof of Theorem \ref{theorem2.8}:}\par
We begin with the case $q_1=+\infty$. Let $0<\varepsilon<1$. If $L_s^kf(z)=(1-|z|^2)^{k-s} (I+R)^kf(z)$, the fact that   $(I+R)^kf$ is holomorphic gives that
$$|L_s^kf(r\zeta)|\preceq \left( \frac1{(1-r^2)^{n+1}}\int_{K(r\zeta,c(1-r^2))} |(I+R)^kf(z)|^\varepsilon dv(z)\right)^\frac1{\varepsilon} (1-r^2)^{k-s},$$
where for $y\in\B^n$ $K(y,t)$ is the nonisotropic ball in $\B^n$ given by
$$K(y,t)=\{ z\in\B^n\,;\, |\overline{z}(z-y)|+|\overline{y}(y-z)|<t\,\}.$$
In \cite{ortegafabrega} it is obtained that
$$|L_s^kf(r\zeta)|\preceq\left( M_{HL}\left( \int_0^1 |(I+R)^k f(t\eta)|^q (1-t^2)^{(k-s)q-1} dt\right)^\frac{\varepsilon}{q}(\zeta)\right)^\frac1{ \varepsilon}.$$
Thus 
\begin{equation*}\begin{split}
&||f||_{HF_s^{p\infty}(w)}^p=\int_{{\bf S}^n} \sup_{0<r<1} |L_s^kf(r\zeta)|^p w(\zeta)d\sigma(\zeta)\\
& \preceq\int_{{\bf S}^n} \left( M_{HL}\left( \int_0^1 |(I+R)^k f(t\eta)|^q (1-t^2)^{(k-s)q-1} dt\right)^\frac{\varepsilon}{q}(\zeta)\right)^\frac{p}{ \varepsilon}w(\zeta)d\sigma(\zeta).
\end{split}\end{equation*}
Since $\frac{p}{\varepsilon}>p$, and $w$ is an $A_p$-weight, $w$ is in $A_{\frac{p}{\varepsilon}}$, and in consequence the unweighted Hardy-Littlewood maximal function is a bounded map $L^\frac{p}{\varepsilon}(w)$ to itself. Hence the above is bounded by
$$C\int_{{\bf S}^n}  \left(\int_0^1 |(I+R)^k f(t\zeta)|^q (1-t^2)^{(k-s)q-1} dt \right)^\frac{p}{q}w(\zeta)d\sigma(\zeta)= C||f||_{HF_s^{ pq}(w)}^p.$$

In order to finish the theorem, we will prove that if $q_0<q_1<+\infty$,
then
$$||f||_{HF_s^{ pq_1}(w)}\leq ||f||_{HF_s^{ pq_0}(w)}^\frac{q_0}{q_1} ||f||_{HF_s^{ p\infty}(w)}^{1-\frac{q_0}{q_1}}.$$
Since 
\begin{equation*}\begin{split}
&||f||_{HF_s^{ pq_1}(w)}^p\leq\int_{{\bf S}^n} \left(\sup_{0<r<1} |(I+R)^kf(r\zeta)| (1-r)^{k-s}\right)^{(q_1-q_0)p/q_1}\\
&\times\left( \int_0^1 |(I+R)^k f(r\zeta)|^{q_0} (1-r^2)^{(k-s)q_0-1}dr\right)^\frac{p}{q_1} w(\zeta)d\sigma(\zeta),
\end{split}\end{equation*}
H\"older's inequality with exponent $q_1/q_0>1$, gives that the above is bounded by
$$C ||f||_{HF_s^{ pq_0}(w)}^{p\frac{q_0}{q_1}}||f||_{HF_s^{ p\infty}(w)}^{p(1- \frac{q_0}{q_1})}.\qed
$$

We now consider  the weighted Hardy space $H^p(w)$, for $1<p<+\infty$, and $w$ an $A_p$ weight.
It is shown in \cite{luecking} that $f\in H^p(w)$ if and only if $f=C[f^*]$, where   $f^*(\zeta)=\lim_{r\rightarrow1}f(r\zeta)\in L^p(w)$ is the radial limit, $C$ is the Cauchy-Szeg\"o kernel.  In addition,  $f=P[f^*]$, where $P$ is the Poisson-Szeg\"o kernel. It follows also that $
||f||_{H^p(w)}^p\simeq ||f^*||_{L^p(w)}$.

 It is immediate to deduce from this that $f\in H^p(w)$ if and only if for any $\alpha \geq 1$, $M_\alpha (f)\in L^p(w)$,
   where $M_\alpha$ is the $\alpha$-admissible maximal operator given by
 $$M_\alpha(f)(\zeta)=\sup_{z\in D_\alpha(\zeta)}|f(z)|.$$
  In addition 
 $||f||_{H^p(w)}\simeq ||M_\alpha(f)||_{L^p(w)}$, with constant that   depends on $\alpha$. Indeed, since $|f(r\zeta)|\leq M_\alpha(f)(\zeta)$, we have that $||f||_{H^p(w)}\leq ||M_\alpha(f)||_{L^p(w)}$. On the other hand, assume that $f\in H^p(w)$. Then $f=P[f^*]$, $f^*\in L^p(w)$ and since $M_\alpha(f)\leq CM_{HL}(f^*)$, (see for instance \cite{rudin}), we deduce that
 \begin{equation*}\begin{split}&\int_{{\bf S}^n} (M_\alpha(f)(\zeta))^p w(\zeta)d\sigma(\zeta) \\&\preceq
 \int_{{\bf S}^n} (M_{HL}(f^*)(\zeta))^p w(\zeta)d\sigma(\zeta)\preceq \int_{{\bf S}^n} |f^* (\zeta)|^p w(\zeta)d\sigma(\zeta)\preceq 
   ||f||_{H^p(w)}^p,
 \end{split}\end{equation*}
where we have used that since $w$ in an $A_p$-weight, the Hardy-Littlewood maximal operator maps $L^p(w)$ continuously  to itself.

Our next result gives a proof for the weighted nonisotropic case of the fact that the spaces $H^p(w)$ can also be defined in terms of admissible area functions.  Similar results, but using a different approach  based on localized good-lambda inequalities, have been obtained in \cite{strombergtorchinsky}  for weighted isotropic Hardy spaces in $\R^n$ .
\begin{theorem}\label{theorem2.9}
Let $1<p<+\infty$, and $w$ be an  $A_p$-weight. Let $f$ be an holomorphic function on $\B^n$. Then the following assertions are equivalent:

(i) $f$ is in $H^p(w)$.

(ii) There exists $\alpha\geq 1$, $k>0$, such that $A_{\alpha, k,2,0}(f)\in L^p(w)$.

 (iii) For every $\alpha\geq 1$,  and $k>0$,   $A_{\alpha, k,2,0}(f)\in L^p(w)$.

In addition, there exists $C>0$ such that for any $f\in H^p(w)$,

$$\frac1{C}||f||_{H^p(w)}\leq ||A_{\alpha, 1,2,0}(f)||_{L^p(w)}\leq C ||f||_{H^p(w)}.$$
\end{theorem} 
{\bf Proof of Theorem \ref{theorem2.9}:}\par

 We already know that (ii) and (iii) are equivalent, so we only have to check the equivalence of (i) and (ii) for the case $k=1$. The proof of (i) implies (ii) is given in \cite{kangkoo}, using the arguments of \cite{stein2}. For the proof of (ii) implies (i), we will follow some ideas of \cite{ahernbrunacascante}.
 
 Without loss of generality we may assume that $f(0)=0$. Let us assume first that $f\in H(\overline{\B^n})$. Then $f=P[f^*]$ where $f^*\in {\mathcal C}({\bf S}^n)$. We want to check that 
 $$||f^*||_{L^p(w)}\leq C||A_{\alpha,1,2,0}(f) ||_{L^p(w)}.$$
 
 We will use that the dual space of $L^p(w)$ can be identified with $L^{p'}(w^{-(p'-1)})$ if the duality is given by
 $$<f,g>=\int_{{\bf S}^n}f(\zeta)\overline{g(\zeta)}d\sigma(\zeta).$$

 Hence, 
 $$||f^*||_{L^p(w)} =\sup\{ |\int_{{\bf S}^n} f^*(\zeta) g^*(\zeta)  d\sigma(\zeta)|,\, g^*\in {\mathcal C}({\bf S}^n),\, ||g^*||_{L^{p'}(w^{-(p'-1)})}\leq 1\,\}.$$
If $g=P[g^*]$, we have (see \cite{ahernbrunacascante} page 131)
\begin{equation}\begin{split}\label{equationformula}
&\frac{n\pi^n}{(n-1)!}\int_{{\bf S}^n} f^*(\zeta)g^*(\zeta)d\sigma(\zeta)\\
&=n^2\int_{\B^n} f(z)g(z)dv(z)+\int_{\B^n} (\nabla_{\B^n}f(z),\nabla_{\B^n}g(z))_{\B^n} \frac{dv(z)}{1-|z|^2},
\end{split}\end{equation}
 where $\nabla_{\B^n}$ is the  gradient in the Bergman metric (see for instance \cite{stein2}), and
 $$(F(z),G(z))_{\B^n}= (1-|z|^2)(\sum_{i,j} (\delta_{i,j}-z_i\overline{z_j}) F_i(z)\overline{G}_j(z)).$$
We then have (see \cite{stein2}) that since $F$ is holomorphic
 \begin{equation*}\begin{split}&||\nabla_{\B^n}F(z)||_{\B^n}^2= (\nabla_{\B^n}F(z), \nabla_{\B^n}F(z))_{\B^n}\\& \simeq (1-|z|^2)\left\{ \sum_{i=1}^n|\frac{\partial}{\partial z_i}F(z)|^2- |\sum_{i=1}^n z_i \frac{\partial}{\partial z_i}F(z)|^2\right\}.\end{split}\end{equation*}
   In order to estimate $\int_{\B^n}f(z)g(z)dv(z)$ we will need to obtain estimates of the values of the functions $f,g$ on compact subsets of $\B^n$ in terms of  the norms $||A_{\alpha, 1,2,0}(f)||_{L^p(w)}$ and $||A_{\alpha, 1,2,0}(g)||_{L^{p'}(w^{-(p'-1)})}$ respectively.
   
   \begin{lemma}\label{lemma2.10}
   Let $1<p<+\infty$ and $w$ an $A_p$-weight. There exists $C>0$ such that for any holomorphic function $f$ in $\B^n$, and any $z=r\zeta$ 
   $$|f(z)|\preceq \left( |f(0)|+\int_0^r \frac{dt}{W(B(\zeta, 1-t^2))^\frac1{p}(1-t^2)}||A_{\alpha, 1,2,0}(f)||_{L^p(w)}^p \right).$$
   
   In particular, if $K\subset \B^n$ is compact and $$||f||_K=\sup_{z\in K}|f(z)|,$$ then there exists a constant $C>0$, depending only on $w$, $p$ and $K$ such that $||f||_K\leq C\left( |f(0)|+ ||A_{\alpha, 1,2,0}(f)||_{L^p(w)}\right)$.
   \end{lemma}
   {\bf Proof of Lemma \ref{lemma2.10}:}\par
   Since $f$ is holomorphic,  we obtain that
   if $z=r\zeta\in\B^n$, there exist $C_i>0$, $i=1,2$, such that for any $\eta\in B(\zeta, C_1(1-r^2))$, then
   \begin{equation*}\begin{split}
   &|\nabla f(z)|^2\preceq \frac{1}{(1-|z|^2)^{n+1}}\int_{K(z,C_2(1-|z|^2))} |\nabla f(y)|^2dv(y) \\
   & \preceq\frac{1}{(1-|z|^2)^2}\int_{K(z,C_2(1-|z|^2))} (1-|y|^2)^{1-n}|\nabla f(y)|^2dv(y) \leq \frac{C}{(1-|z|^2)^2} (A_{\alpha, 1,2,0}(f)(\eta))^2.
   \end{split}\end{equation*}
   Consequently
   $$\left((1-|z|^2)|\nabla f(z)|\right)^p\preceq  ( A_{\alpha, 1,2,0}(f)(\eta))^p.$$
  Then we have
    \begin{equation*}\begin{split}&\left((1-|z|^2)|\nabla f(z)|\right)^pW(B(\zeta, 1-r^2))\preceq  
   \int_{B(\zeta, 1-r^2)} ( A_{\alpha, 1,2,0}(f)(\eta))^pw(\eta)d\sigma(\eta) \preceq\\& || A_{\alpha, 1,2,0}(f)||_{L^p(w)}^p.\end{split}\end{equation*}
   In particular, if $0<r<1$ and $\zeta\in{\bf S}^n$,
   $$|\frac{\partial f}{\partial r}(r\zeta)|\preceq   \frac{1}{W(B(\zeta, 1-r^2))^\frac1{p}(1-r^2)}||A_{\alpha, 1,2,0}(f)||_{L^p(w)},
   $$
   and integrating, we finally obtain
   $$|f(r\zeta)|\preceq \left( |f(0)|+\int_0^r\frac{dt}{W(B(\zeta, 1-t^2))^\frac1{p}(1-t^2)}||A_{\alpha, 1,2,0}(f)||_{L^p(w)}\right). $$

   For the remaining affirmation, let $K\subset\B^n$ be compact. Then there exists $0<\delta<1$ such that for any $z=r\zeta\in K$, $r\leq 1-\delta$, and  
    $$|f(z)|\preceq \left( |f(0)|+ \frac{1}{W(B(\zeta, \delta))^\frac1{p}\delta }||A_{\alpha, 1,2,0}(f)||_{L^p(w)}\right).$$
    Since $w$ is doubling, and there exists $N>0$ (not depending on $\zeta $) such that ${\bf S}^n\subset B(\zeta, cN\delta))$, $W({\bf S}^n)\preceq  W(B(\zeta, \delta))$, and consequently $$||f||_K\preceq    |f(0)|+ ||A_{\alpha, 1,2,0}(f)||_{L^p(w)} .\qed$$
    
    Going back to the proof of the Theorem \ref{theorem2.9},  let $0<\varepsilon<1$. The above lemma together with the fact that if $w$ is an $A_p$ weight, then $w^{-(p'-1)}$ is an $A_{p'}$-weight, give  by (\ref{equationformula}) that
 \begin{equation*}\begin{split}
  & | \int_{{\bf S}^n} f^*(\zeta)g^*(\zeta)d\sigma(\zeta)| \\
& \preceq||A_{\alpha, 1,2,0}(f)||_{L^p(w)}||A_{\alpha, 1,2,0}(g)||_{L^{p'}(w^{-(p'-1)})} +|\int_{1-\varepsilon\leq |z|<1} f(z)g(z)dv(z)|\\
&+\int_{\B^n} ||\nabla_{\B^n}f(z)||_{\B^n}||\nabla_{\B^n}g(z)||_{\B^n} \frac{dv(z)}{1-|z|^2}.
 \end{split}\end{equation*}
    
   In order to estimate the second integral, we use polar coordinates, and obtain 
  $$|\int_{1-\varepsilon\leq |z|<1} f(z)g(z)dv(z)|,$$
  which by H\"older's inequality is bounded by 
    \begin{equation*}\begin{split}
   &\int_{1-\varepsilon}^1\int_{{\bf S}^n} |f(r\zeta)||g(r\zeta)|d\sigma(\zeta) dr \\
   & \preceq\int_{1-\varepsilon}^1 ||f_r||_{L^p(w)}||g_r||_{L^{p'}(w^{-(p'-1)})} dr\preceq \varepsilon ||f||_{H^p(w)}||g||_ {H^{p'}(w^{-(p'-1)})}\\
   &\preceq   
    \varepsilon ||f^*||_{L^p(w)}||g^*||_ {L^{p'}(w^{-(p'-1)})}.
   \end{split}\end{equation*}
   For the third integral, we use (5.1) of \cite{coifmanmeyerstein} to estimate it by
   $$\int_{{\bf S}^n}A_{\alpha, 1,2,0}(f)(\zeta) A_{\alpha, 1,2,0}(g)(\zeta)d\sigma(\zeta) \preceq ||A_{\alpha, 1,2,0}(f)||_{L^p(w)} ||A_{\alpha, 1,2,0}(g)||_{L^{p'}(w^{-(p'-1)})}.$$
   
   Since we already know (see \cite{kangkoo}) that $||A_{\alpha, 1,2,0}(g)||_{L^{p'}(w^{-(p'-1)})}\preceq   ||g^*||_{L^{p'}(w^{-(p'-1)})}$, we finally obtain
   $$||f^*||_{L^p(w)}\preceq  ||A_{\alpha, 1,2,0}(f)||_{L^p(w)}+  \varepsilon ||f^*||_{L^p(w)},$$
   which gives the result for $f\in H(\overline{\B}^n)$.
   
   So we are left to show that the estimate we have already obtained holds for a general holomorphic function in $\B^n$.
   If $f$ is an holomorphic function on $\B^n$ such that $||A_{\alpha,1,2,0}(f)||_{L^p(w)}<+\infty$, let $f_r(z)=f(rz)\in H(\overline{\B^n})$, for $0<r<1$. We then have that
   \begin{equation}\label{areaformula}
   ||f_r||_{H^p(w)}\preceq ||A_{\alpha, 1,2,0}(f_r)||_{L^p(w)}
   \end{equation}
   Let us check first that $$\sup_{r}||A_{\alpha, 1,2,0}(f_r)||_{L^p(w)}\leq C ||A_{\alpha, 1,2,0}(f)||_{L^p(w)}.$$
   
     Notice that
  $$||A_{\alpha, 1,2,0}(f_r)||_{L^p(w)}^p= ||J_\alpha((1-|\cdot|^2) (I+R)f_r)||_{ L^p(w)(L^2(\frac{dv(z)}{(1-|z|2)^{n+1}}))}.$$  
  
  We will check that there exists $0\leq G(\zeta,z)\in   L^p(w)(L^2(\frac{dv(z)}{(1-|z|^2)^{n+1}}))$ such that for any $0<r<1$, $\zeta\in{\bf S}^n$, $z\in\B^n$, $J_\alpha ((1-|\cdot|^2)(I+R)f_r)(\zeta, z)\leq  G(\zeta,z)$, and $||G||_{ L^p(w)(L^2(\frac{dv(z)}{(1-|z|^2)^{n+1}}))}\preceq ||A_{\alpha,1,2,0}(f)||_{L^p(w)}$. 
  
  Let us obtain such a function $G$. Since   by hypothesis $A_{\alpha,1,2,0}f\in  L^p(w)$, we have  that the holomorphic function $f$ satisfies that  $A_{\alpha,1,2,0}f\in L^1(d\sigma)$, and consequently that there exists $C>0$ such that for any $z\in \B^n$, $|f(z)|\preceq  \frac{1}{(1-|z|^2)^n}$. Hence, the integral representation theorem gives that for $N>0$ sufficiently large, and $z\in\B^n$,
  $$(I+R)f(rz)=C\int_{\B^n} \frac{(1-|y|^2)^N (I+R)f(y)}{(1-rz\overline{y})^{n+1+N}}dv(y).$$
  Next, there is a constant $C>0$ such that for any $0<r<1$, $z,y\in\B^n$,
  $|1- rz\overline{y}|\geq C |1- z\overline{y}|$, and the above formula gives that
  $$|(I+R)f(rz)|\preceq \int_{\B^n} \frac{(1-|y|^2)^N |(I+R)f(y)|}{|1- z\overline{y}|^{n+1+N}}dv(y).$$
  Combining the above results we have that 
  \begin{equation*}\begin{split}
  &\chi_{D_\alpha(\zeta)}(z)(1-|z|^2)|(I+R)f(rz)| \\
  &  \preceq\chi_{D_\alpha(\zeta)}(z)\int_{\B^n} \frac{(1-|y|^2)^{N-1}(1-|z|^2) ((1-|y|^2)|(I+R)f(y)|)}{|1- z\overline{y}|^{n+1+N}}dv(y)\\
  &= C\chi_{D_\alpha(\zeta)}(z) P^{N-1,1}((1-|\cdot|^2)(I+R)f)(z):=G(z,\zeta).
  \end{split}\end{equation*}
  Theorem \ref{theorem2.5} shows that provided $N$ is chosen sufficiently large, $P^{N-1,1}$ maps $F^{\alpha,p,2}(w)$ to itself, and in particular that
  \begin{equation*}\begin{split}&||G||_{ L^p(w)(L^2(\frac{dv(z)}{(1-|z|^2)^{n+1}}))} = ||P^{N-1,1}((1-|\cdot|^2)(I+R)f)||_{\alpha, p,2,w} \\
  &\preceq || (1-|\cdot|^2)(I+R)f||_{\alpha, p,2,w}= C ||A_{\alpha,1,2,0}(f)||_{L^p(w)}<+\infty.
  \end{split}\end{equation*}
  
   Consequently $$||f_r||_{H^p(w)} \preceq  ||A_{\alpha, 1,2,0}(f )||_{ L^p(w)},$$  and therefore $f\in H^p(w)$.
    \qed

We will now remark on some facts about weighted Hardy-Sobolev spaces. Let us recall, that if $1<p<+\infty$, $0<s<n$, and $w$ is an $A_p$-weight, we denote by $H_s^p(w)$ the space of holomorphic functions $f$ on $\B^n$ satisfying that    
   $$||f||_{H_s^p(w)} =||(I+R)^sf||_{H^p(w)}<+\infty.$$ The results obtained in the previous theorems give alternative equivalent definitions of the spaces $H_s^p(w)$ in terms of admissible maximal or radial functions and admissible area functions. 
   
   On the other hand, when $w\equiv1$, and $0<s<n$, it is well known, see for instance \cite{cascanteortegat}, that the space $H_s^p$ admits a representation in terms of a   fractional Cauchy-type kernel $C_s$ defined by
   $$C_s(z,\zeta)=\frac1{(1-z\overline{\zeta})^{n-s}}.$$ 
  
   The same lines of the proof of the unweighted case can be used  to obtain  a similar characterization in the weighted case. We just have to use that the Hardy-Littlewood maximal operator is bounded in $ L^p(w)$, if $w$ is an $A_p$-weight and Lemma \ref{weightedspaces}.
   \begin{theorem}\label{characterizationhardysobolev}
    Let $1<p<+\infty$, $0<s<n$, and $w$ be an $A_p$-weight. We then have that the map $$C_s(f)(z)=\int_{\B^n} \frac{f(\zeta)}{(1-z\overline{\zeta})^{n-s}}d\sigma(\zeta),$$
     is a bounded map of $ L^p(w)$ onto $H_s^p(w)$.
    
   \end{theorem} 
    
\section{Holomorphic potentials and Carleson measures}

In this section we will study Carleson measures for $H_s^p(w)$, $1<p<+\infty$ and $0<s<n$, that is, the  positive finite Borel measures $\mu$ on $\B^n$ satisfying
\begin{equation}\label{equation3.1}
\| f\|_{L^p(d\mu)}\leq C  \, \|f\|_{H_s^p(w)}, \qquad f \in
 L^p(w).
\end{equation} 
  In what follows we will write $${\displaystyle{\int\!\!\!\!\!\setminus}_Ewd\sigma=\frac{1}{|E|}\int_Ew},$$ where $E$ is a measurable set in ${\bf S}^n$, $|E|$ denotes its Lebesgue measure.
  
  By Theorem \ref{characterizationhardysobolev}, this inequality can be rewritten as follows:
  \begin{equation}\label{equation3.2}
\| C_s( f)\|_{L^p(d\mu)}\leq C \, \|f\|_{ L^p(w)}, \qquad f \in
 L^p(w).  
\end{equation}  
 We recall that we have defined the non-isotropic potential of a positive Borel function $f$ on ${\bf S}^n$ by
 \begin{equation}\label{equation3.3}
 K_s(f)(z)=\int_{{\bf S}^n}K_s(z,\zeta)f(\zeta)d\sigma(\zeta)=\int_{{\bf S}^n}\frac{f(\zeta)}{|1-z\overline{\zeta}|^{n-s}}d\sigma( \zeta),
 \end{equation} 
 for $z\in\overline{\B}^n$. 
 
 Analogously to what happens for isotropic potentials (see \cite{adams}), in the nonisotropic case it can be proved that if $w$ is an $A_p$ weight and $\zeta_0\in{\bf S}^n$ satisfies that \begin{equation}\label{continuity}\int_{{\bf S}^n}\frac1{|1-\zeta_0\overline{\zeta}|^{(n-s)p'}}w^{-(p'-1)}(\zeta) d\sigma(\zeta)<+\infty,\end{equation}  then for any $f\in  L^p(w)$, $K_s(f)$ is continuous in $\zeta_0$.  
 Observe that when $w\equiv 1$, (\ref{continuity}) holds if and only if $n-sp<0$. In the general weighted case, if $w$ satisfies a doubling condition of order $\tau$,  and $\tau-sp<0$, we also have that (\ref{continuity}) holds, and consequently the Carleson measures in this case for weighted Hardy Sobolev spaces are just the finite ones. 
 Indeed, assume that $\tau-sp<0$. We then have
 
\begin{equation*}\begin{split}
& \int_{{\bf S}^n}\frac1{|1-\zeta_0\overline{\zeta}|^{(n-s)p'}}w^{-(p'-1)}(\zeta) d\sigma(\zeta)=\int_{{\bf S}^n} w^{-(p'-1)}(\zeta) \int_{|1-\zeta_0\overline{\zeta}|<t}\frac{dt}{t^{(n-s)p'}}d\sigma(\zeta)\\
&\leq\int_0^K \frac{\int_{B(\zeta_0,t)}  w^{-(p'-1)}}{t^{(n-s)p'}}\simeq \int_0^K \frac{t^ndt}{\left({\int\!\!\!\!\!\setminus}_{B(\zeta_0,t)} w \right)^{p'-1}t^{(n-s)p'}}\preceq \sum_k \frac{2^{-ksp'}}{W(B(\zeta_0, 2^{-k}))}.
\end{split}\end{equation*}
The fact that $w$ satisfies condition $D_\tau$ gives that $W({\bf S}^n)\preceq 2^{k\tau} W(B(\zeta_0, 2^{-k}))$, and consequently the above sum is bounded, up to constants, by
$$\sum_k  2^{k(\tau(p'-1)-sp'} .$$
 Since $\tau-sp<0$ we also have that $\tau(p'-1)-sp'<0$, and we are done.
 
 From now on we will assume that $\tau-sp\geq0$.

 The problem of characterizing the positive finite Borel measures $\mu$ on $\B^n$ for which the following  inequality holds
 \begin{equation}\label{equation3.4}
 || K_s (f)||_{L^p(d\mu)}\leq C\, ||f||_{ L^p(w)},
\end{equation}
has been thoroughly studied, and there are, among others, characterizations in terms of weighted nonisotropic Riesz capacities that are  defined as follows:
if $E\subset{{\bf S}}^n$, $1<p<+\infty$ and $0<s<n$, 
$$C_{sp}^w(E)= \inf \{||f||_{ L^p(w)}^p\,;\, f\geq0,\, K_s(f)\geq 1 \,\,{\rm on}\,\, E\,\}.$$
It is well known, that when $w\equiv 1$, $C_{sp}(B(\zeta,r))\simeq r^{n-sp}$, $\zeta\in{{\bf S}^n}$, $r<1$. See \cite{adams} for expressions of weighted capacities of balls in $\R^n$.

As it happens in $\R^n$ (see \cite{adams}), we have  that if $0\leq n-sp$, (\ref{equation3.4}) holds if and only if there exists $C>0$ such that for any open set $G\subset {\bf S}^n$, 
\begin{equation}\label{equation3.5}
\mu(T(G))\leq C C_{sp}^w(G).
\end{equation}
Here $T(G)$ is the admissible tent over $G$, defined by 
$$T(G)=T_\alpha (G)=\left( \bigcup_{\zeta\notin G} D_\alpha(\zeta)\right)^c.$$

The problem of characterizing the Carleson measures $\mu$ for   the holomorphic case (\ref{equation3.2}) is much more complicated, even in the nonweighted case.  
Since $|C_s(z,\zeta)|\leq K_s(z,\zeta)$, it follows from Theorem \ref{characterizationhardysobolev}, that (\ref{equation3.4}) implies (\ref{equation3.2}), and consequently that if condition (\ref{equation3.5}) is satisfied, then $\mu$ is a Carleson measure for $H_s^p(w)$. Of course, when $n-s<1$ both problems are equivalent, even in the weighted case, simply because if $f\geq 0$, $|C_s(f)|\simeq K_s(f)$,  but  when $n>1$   (see \cite{ahern} and \cite{cascanteortega2}), condition (\ref{equation3.4}) for the unweighted case is not, in general, equivalent to condition (\ref{equation3.2}).  Observe that when $n-sp\leq 0$, $H_s^p$ consists of regular functions, and consequently any finite measure is a Carleson measure for the holomorphic and the real case. It is proved in \cite{cohnverbitsky2} that this equivalence still remains true if we are not too far from the regular case, namely, if  $0\leq n-sp<1$ . The main purpose of this section is to obtain a  result in this line for a wide class of $A_p$-weights.

 In \cite{ahern} it is also shown that if (\ref{equation3.2}) holds for $w\equiv 1$, then the capacity condition on balls is satisfied, i.e. there exists $C>0$ such that $\mu(T(B(\zeta,r)))\leq C r^{n-sp}$, for any $\zeta\in{{\bf S}}^n$ and any $0<r<1$. The   following proposition   obtains a necessary condition in this line for the weighted holomorphic trace inequality.
 \begin{proposition}\label{proposition3.1}
 Let $1<p<+\infty$, $0<s<n$ . Let $\mu$ be a positive finite Borel measure on $\B^n$, and $w$ be an $A_p$-weight.  Assume that there exists $C>0$ such that  $$||f||_{L^p(d\mu)}\leq C||f||_{H_s^p(w)},
 $$
 for any $f\in H_s^p(w)$. We then have that there exists $C>0$ such that for any $\zeta\in{{\bf S}}^n$, $r>0$,
 $$\mu(T(B(\zeta,r))\leq C \frac{W(B(\zeta,r))}{r^{sp}}.$$
 \end{proposition}
 {\bf Proof of Proposition \ref{proposition3.1}:}\par
Let $\zeta\in{{\bf S}}^n$, $0<r<1$ be fixed. If $z\in \overline {\B}^n$, let 
$$F(z)=\frac{1}{(1-(1-r)z\overline{\zeta})^N},$$
with $N>0$ to be chosen later.
If $z\in T(B(\zeta,r))$, and $z_0=\frac{z}{|z|}$, $(1-|z|)\preceq r$ and $|1-z_0\overline{\zeta}|\preceq r$. Hence 
$ 
|1-(1-r)z\overline{\zeta}|\preceq r,
$ 
and consequently, 
$$
\frac{\mu(T(B(\zeta,r)))}{r^{Np}}\leq C \int_{T(B(\zeta,r))} |F(z)|^pd\mu(z).
$$
On the other hand, 
\begin{equation*}\begin{split}
&||F||_{H_s^p(w)}^p\leq C\int_{{\bf S}^n}\frac1{|1-(1-r)\eta\overline{\zeta}|^{(N+s)p}} w(\eta)d\sigma(\eta)\\
&=\int_{B(\zeta ,r)}\frac1{|1-(1-r)\eta\overline{\zeta}|^{(N+s)p}}w( \eta)d\sigma(\eta)\\&+\sum_{k\geq 1} \int_{B(\zeta,2^{k+1}r)\setminus B(\zeta,2^{k}r) }\frac1{|1-(1-r)\eta\overline{\zeta}|^{(N+s)p}}w( \eta)d\sigma(\eta).
\end{split}\end{equation*}
If $k\geq 1$, and $\eta\in B(\zeta,2^{k+1}r)\setminus B(\zeta,2^{k}r) $, $|1-(1-r)\eta\overline{\zeta}|\simeq 2^kr$. This estimates together with the fact that   $w$ is doubling, give that the above is bounded by
 \begin{equation*}\begin{split}
 &\sum_{k\geq 0}\frac{W(B(\zeta,2^{k+1}r))}{(2^k r)^{(N+s)p}}\preceq \frac{W(B(\zeta,r))}{r^{(N+s)p}}\sum_{k\geq 0}\left( \frac{C}{2^{(N+s)p}}\right)^k,
 \end{split}\end{equation*}
 which gives the desired estimate, provided $N$ is chosen big enough.\qed
 
 We observe that for some special weights besides the case $w\equiv1$,  the expression that appears in the above proposition $\frac{W(B(\zeta,r))}{r^{sp}}$ coincide with the weighted capacity of a ball  (see \cite{adams}).

 If $\nu$ is a positive Borel measure on ${{\bf S}}^n$,   $1<p<+\infty$, $0<s<n$ and $w$ is an $A_p$-weight, it is introduced in \cite{adams} the $(s,p)$-energy of $\nu$ with weight $w$, which is defined by
 \begin{equation}\label{weightedenergy}
 {\mathcal E}_{sp}^w(\nu)=\int_{{\bf S}^n} (K_s(\nu)(\zeta))^{p'}w(\zeta)^{-(p'-1)}d\sigma(\zeta).
 \end{equation} 
 If we write $(K_s(\nu))^{p'}= (K_s(\nu))^{p'-1} K_s(\nu)$, Fubini's theorem gives that
 $${\mathcal E}_{sp}^w(\nu)=\int_{{\bf S}^n} {\mathcal U}_{sp}^w(\nu)(\zeta) d\nu(\zeta),
 $$
 where
 $${\mathcal U}_{sp}^w(\zeta)=K_s(w^{-1}  K_s(\nu))^{p'-1}(\zeta)$$ is the weighted nonlinear potential of the measure $\nu$.
 When $w\equiv 1$, Wolff's theorem (see \cite{hedbergwolff}) gives another representation of the energy, in terms of the so-called Wolff's potential.
 
 In the general case, it is introduced in \cite{adams} a weighted Wolff-type potential of a measure $\nu$ as
 \begin{equation}\label{weightedwolffpotential}
 {\mathcal W}_{sp}^w(\nu)(\zeta)= \int_0^1\left( \frac{\nu(B(\zeta,1-r))}{(1-r)^{n- sp}}\right)^{p'-1} {\int\!\!\!\!\!\setminus}_{B(\zeta, 1-r)}w^{-(p'-1)}(\eta)d\sigma(\eta)\frac{dr}{1-r}.
 \end{equation}
 In the same paper, it is shown that provided $w$ is an $A_p$-weight, the following weighted Wolff-type theorem holds:
 \begin{equation}\label{weightedwolfftheorem}
 {\mathcal E}_{sp}^w(\nu)\simeq \int_{{\bf S}^n} {\mathcal W}_{sp}^w(\nu)(\zeta) d\nu(\zeta).
 \end{equation}
 In fact, we have the pointwise estimate ${\mathcal W}_{sp}^w(\nu)(\zeta)\leq C{\mathcal U}_{sp}^w(\nu)(\zeta)$, and Wolff's theorem gives that the converse is true, provided we integrate with respect to $\nu$.

 In \cite{adams}  a weighted extremal theorem for the weighted Riesz capacities it is also shown, namely, if $G\subset{{\bf S}}^n$ is open, there exists a positive capacitary measure $\nu_G$ such that
 
 (i) ${\rm supp}\,\, \nu_G\subset G$.
 
 (ii) $\nu_G(G)=C_{sp}^w(G)={\mathcal E}_{sp}^w(\nu_G)$.
 
 (iii) ${\mathcal W}_{sp}^w(\nu_G)(\zeta)\geq C$, for $C_{sp}^w$-a.e. $\zeta\in G$.
 
 (iv) ${\mathcal W}_{sp}^w(\nu_G)(\zeta)\leq C$, for any $\zeta\in {\rm supp}\,\,\nu_G$.
 
 We now introduce two holomorphic weighted Wolff-type potentials, which generalize the ones defined in \cite{cohnverbitsky2}. These potentials will be used in the proof of the characterization of the Carleson measures for $H_s^p(w)$, for the case $0\leq\tau-sp<1$. 
 Let $1<p<+\infty$, $0<s<\frac{n}{p}$, and $\nu$ be a positive Borel measure on ${{\bf S}}^n$. For any $\lambda>0$, and $z\in\B^n$, we set
 \begin{equation}\begin{split}\label{holompot1}
 &{\mathcal U}_{sp}^{w\lambda}(\nu)(z)\\
 &=\int_0^1 \int_{{\bf S}^n}\left( \frac{\nu(B(\zeta, 1-r))}{(1-r)^{n-sp}}\right)^{p'-1} \frac{(1-r)^{\lambda-n}}{ (1-rz\overline{\zeta})^\lambda} \left({\int\!\!\!\!\!\setminus}_{B(\zeta, 1-r)}w^{-(p'-1)}\right) d\sigma(\zeta) \frac{dr}{1-r},
 \end{split}\end{equation}
 and
 \begin{equation}\begin{split}\label{holompot2}
 &{\mathcal V}_{sp}^{w\lambda}(\nu)(z)\\
 &=\int_0^1 \left(\int_{{\bf S}^n} \frac{(1-r)^{\lambda+sp-n}}{(1-rz\overline{\zeta})^\lambda} \left({\int\!\!\!\!\!\setminus}_{B(\zeta, 1-r)}w^{-(p'-1)}\right)^\frac1{p'-1} d\nu(\zeta) \right)^{p'-1}  \frac{dr}{1-r}.
 \end{split}\end{equation}
 
 Obviously, both potentials are holomorphic functions in the unit ball. We will see, that if $p\leq 2$ the first one is bounded from below by the weighted Wolff-type potential we have just introduced, whereas if $p\geq 2$, the second one is bounded from below by the same potential.
 
 In the unweighted case, \cite{cohnverbitsky2} the proof of the estimates of the holomorphic potentials, rely on an extension of Wolff's theorem. This extension gives that
 if $1<p<+\infty$,  $s>0$, $0<q<+\infty$, and $\nu$ is a positive Borel measure on ${{\bf S}}^n$, then
 $$
 \int_{{\bf S}^n} \left( \int_0^1 \left( \frac{\nu(B(\zeta,t))}{t^{n-s}} \right)^q  \frac{dt}{t}\right)^\frac{p'}{q}   d\sigma(\zeta)\preceq 
 \int_{{\bf S}^n} {\mathcal W}_{sp}^w(\nu)(\zeta)d\nu(\zeta).
$$
 Observe that if the above estimate holds for one $q_0$, it also holds for any $q\geq q_0$. The case $q=1$ is the integral estimate in Wolff's theorem, since we have that $${\mathcal E}_{sp}(\nu)\simeq \int_{{\bf S}^n}\left( \int_0^1  \frac{\nu(B(\zeta,t))}{t^{n-s}}   \frac{dt}{t}\right)^{p'}   d\sigma(\zeta).$$  
 
 The arguments in \cite{cohnverbitsky2} can easily be used to show the following weighted version   of the above theorem. We omit the details of the proof.
 
 \begin{theorem}\label{extenssionwolfftheorem}
 Let $1<p<+\infty$, $w$ an $A_p$ weight, $s>0$, $K>0$, $0<q<+\infty$, and $\nu$ be a positive Borel measure on ${{\bf S}}^n$. Then
 \begin{equation}\begin{split}\label{formulaextenssionwolfftheorem}
 &\int_{{\bf S}^n} \left( \int_0^K \left( \frac{\nu(B(\zeta,t))}{t^{n-s}} \left( {\int\!\!\!\!\!\setminus}_{B(\zeta, t)}w^{-(p'-1)}(\eta)d\sigma(\eta) \right)^\frac1{p'-1}\right)^q \frac{dt}{t}\right)^\frac{p'}{q} w(\zeta) d\sigma(\zeta) \\
 &\preceq\int_{{\bf S}^n} {\mathcal W}_{sp}^w(\nu)(\zeta)d\nu(\zeta).
 \end{split}\end{equation}
 \end{theorem}

  Before we obtain estimates of the $H_s^p(w)$-norm of the weighted holomorphic potentials already introduced, we will give a characterization  for weights satisfying a doubling condition  
 \begin{lemma}\label{doubling}
 Let $1<p<+\infty$ and $w$ be an $A_p$ weight on ${\bf S}^n$, and assume that $w\in D_\tau$, for some $\tau>0$. We then have:
 
 (i) For any $t\in\R$ satisfying that $t>\tau-n$, there exists $C>0$ such that 
 \begin{equation}\label{diniconditionawayzero}\int_{r}^{+\infty} \frac1{x^t}{\int\!\!\!\!\!\setminus}_{B(\zeta,x)} w \frac{dx}{x}\leq C\frac1{r^t}{\int\!\!\!\!\!\setminus}_{B(\zeta,r)} w ,  \end{equation}
 $r<1$, $\zeta\in{\bf S}^n$.
 
 (ii) For any $t\in\R$ satisfying that $t>\tau-n$, there exists $C>0$ such that 
 \begin{equation} \label{diniconditionnearzero}\int_0^{r}  x^t \left({\int\!\!\!\!\!\setminus}_{B(\zeta,x)} w^{-(p'-1)} \right)^{p-1}\frac{dx}{x}\leq C r^t \left({\int\!\!\!\!\!\setminus}_{B(\zeta,r)} w^{-(p'-1)} \right)^{p-1},  \end{equation}
 $r<1$, $\zeta\in{\bf S}^n$.
 \end{lemma}
  {\bf Proof of Lemma \ref{doubling}:}\par
 We begin with the proof of part (i). Let $t>\tau-n$. Then
 \begin{equation*}\begin{split}
 &\int_r^{+\infty} \frac1{x^{t}}{\int\!\!\!\!\!\setminus}_{B(\zeta,x)} w \frac{dx}{x}= \sum_{k\geq0} \int_{2^kr}^{2^{k+1}r}\frac1{x^{t}}{\int\!\!\!\!\!\setminus}_{B(\zeta,x)}w \frac{dx}{x}  \\
 &\preceq\sum_{k\geq 0} \frac{1}{2^{k(t+n)}r^{t+n}} W(B(\zeta,2^{k+1}r))\preceq   \sum_{k\geq 0} \frac{1}{2^{k(t+n)}r^{t+n}}2^{k\tau}W(B(\zeta,r))= C\frac1{r^{\delta }}{\int\!\!\!\!\!\setminus}_{B(\zeta,r)}w,
\end{split}\end{equation*}
 since $w$ is in $D_\tau$, and $t+n>\tau$.

Next we show that (ii) holds. If $\zeta\in{{\bf S}^n}$ and $r>0$, the fact that $w\in A_p$ gives that $\left({\int\!\!\!\!\!\setminus}_{B(\zeta,x)} w^{-(p'-1)} \right)^{p-1} \simeq  \left({\int\!\!\!\!\!\setminus}_{B(\zeta,x)} w \right)^{-1}$, and consequently,
\begin{equation*}\begin{split}
 &\int_0^r   x^{t} \left({\int\!\!\!\!\!\setminus}_{B(\zeta,x)} w^{-(p'-1)} \right)^{p-1}    \frac{dx}{x}= \sum_{k\geq0} \int_{2^{-k}r}^{2^{-k+1}r} x^{t} \left({\int\!\!\!\!\!\setminus}_{B(\zeta,x)} w^{-(p'-1)} \right)^{p-1}  \frac{dx}{x} \\
 &\preceq \sum_{k\geq 0}  {2^{-kt}r^{t}} \frac1{  {\int\!\!\!\!\!\setminus}_{B(\zeta,2^{-k}r)} w }\preceq   \sum_{k\geq 0} \frac{1}{2^{ k(t+n)}}r^{t-n}2^{k\tau}W(B(\zeta,r))\simeq  r^{t } \left({\int\!\!\!\!\!\setminus}_{B(\zeta,r)} w^{-(p'-1)} \right)^{p-1}.
\end{split}\end{equation*}
  \qed

{\bf Remark:} In fact, it can be proved that both conditions (i) and (ii) are in turn equivalent to the fact that the $A_p$ weight is in $D_\tau$.

 We can now obtain the estimates on the weighted holomorphic potentials defined in (\ref{holompot1}) and (\ref{holompot2}).
  \begin{theorem}\label{pgreaterthan2}
  Let $1< p<+\infty$, $0<\alpha<n$, $w$  an $A_p$-weight. Assume that $w$ is in $D_\tau$ for some $0\leq\tau-sp< 1 $. We then have:
  \begin{enumerate}
   \item If $1<p<2$, there exists $0<\lambda <1$  and $C>0$   such that for any    finite positive Borel measure $\nu$ on ${\bf S}^n$ the following assertions hold:
  \begin{enumerate}
  \item[a)] For any $\eta\in{\bf S}^n$,
  $$\lim_{\rho\rightarrow 1} {\rm Re}\,\, {\mathcal U}_{sp}^{w\lambda}(\nu)(\rho\eta)\geq C {\mathcal W}_{sp}^{w\lambda}(\nu)(\eta).$$
  \item[b)] $|| {\mathcal U}_{sp}^{w\lambda}(\nu)||_{H_s^p(w)}^p \leq C{\mathcal E}_{sp}^w(\nu)$.
  \end{enumerate}
 
  \item If $p\geq 2$, there exists $0<\lambda <1$  and $C>0$   such that for any    finite positive Borel measure $\nu$ on ${\bf S}^n$ the following assertions hold: 
  \begin{enumerate}
  \item[a)] For any $\eta\in{\bf S}^n$,
  $$\lim_{\rho\rightarrow 1} {\rm Re}\,\, {\mathcal V}_{sp}^{w\lambda}(\nu)(\rho\eta)\geq C {\mathcal W}_{sp}^{w\lambda}(\nu)(\eta).$$
  
  \item[b)] $|| {\mathcal V}_{sp}^{w\lambda}(\nu)||_{H_s^p(w)}^p \leq C{\mathcal E}_{sp}^w(\nu)$.
  \end{enumerate}
  \end{enumerate}
  \end{theorem} 
  {\bf Proof of Theorem \ref{pgreaterthan2}:}\par
 We will follow the scheme of \cite{cohnverbitsky2} where it is proved for the unweighted case. The weights introduce   new technical difficulties that require a careful use of the hypothesis $A_p$ and $D_\tau$ that we assume on the weight $w$. In order to make  the proof easier to follow we sketch some of the arguments in \cite{cohnverbitsky2},  emphasizing the necessary changes we need to make in the weighted case.

Let us prove (1).  
We choose $\lambda$ such that $\tau-sp<\lambda<1$ and define ${\mathcal U}_{sp}^{w\lambda}$ as in \ref{holompot1}. Then $\tau-s<\frac{\lambda +s-\tau(2-p)}{p-1}$. Consequently there exists $t$ such that    $\tau-s<t<\frac{\lambda +s-\tau(2-p)}{p-1}$. Observe that   $ t+s-n>\tau-n$ and $\frac{\lambda+s-t(p-1)}{ 2-p}-n>\tau-n$.

We begin now the proof of a). The fact that $\lambda<1$ gives that if $\rho<1$, $\eta\in{\bf S}^n$, and $C>0$,
\begin{equation*}\begin{split}&{\rm Re}\, {\mathcal U}_{sp}^{w\lambda}(\rho\eta) \\&\succeq \int_0^1 \int_{B(\eta, C(1-r))} \left( \frac{\nu(B(\zeta,1-r))}{(1-r)^{n-sp}}\right)^{p'-1} \frac{(1-r)^{\lambda -n}}{|1-r\rho \eta\overline{\zeta}|^{\lambda}}\left({\int\!\!\!\!\!\setminus}_{B(\zeta,1-r)} w^{-(p'-1)}\right)  d\sigma(\zeta)\frac{dr}{1-r}.\end{split}\end{equation*}
If $C>0$ has been chosen small enough, we have that for any $\zeta\in B(\eta,C(1-r))$, $B(\eta,C(1-r))\subset B(\zeta, 1-r)$. In addition, $|1-r\rho\eta\overline{\zeta}|\preceq |1-r\rho|$.  
These estimates, together with the fact that $w^{-(p'-1)}$ satisfies a doubling condition, give that the above integral is bounded from below by
\begin{equation*}\begin{split}
&C\int_0^1 \int_{B(\eta, C(1-r))} \left( \frac{\nu(B(\eta,C(1-r))}{(1-r)^{n-sp}}\right)^{p'-1} \frac{(1-r)^{\lambda -n}}{|1-r\rho |^{\lambda}}\left({\int\!\!\!\!\!\setminus}_{B(\eta,1-r)} w^{-(p'-1)} \right)d\sigma(\zeta)\frac{dr}{1-r}\\
&\geq C\int_{0}^\rho \left( \frac{\nu(B(\eta,C(1-r))}{(1-r)^{n-sp}}\right)^{p'-1} \frac{(1-r)^{\lambda }}{|1-r\rho |^{\lambda}}\left({\int\!\!\!\!\!\setminus}_{B(\eta,1-r)} w^{-(p'-1)}\right) \frac{dr}{1-r} \\&\geq C\int_{0}^\rho \left( \frac{\nu(B(\eta,C(1-r))}{(1-r)^{n-sp}}\right)^{p'-1} \left({\int\!\!\!\!\!\setminus}_{B(\eta,1-r)} w^{-(p'-1)}\right) \frac{dr}{1-r},
\end{split}\end{equation*}
where in last estimate we have used that since $r<\rho$, $1-r\rho\simeq 1-r$.

We have proved then
$$\int_0^\rho \left( \frac{\nu(B(\eta,C(1-r))}{(1-r)^{n-sp}}\right)^{p'-1} \left({\int\!\!\!\!\!\setminus}_{B(\eta,1-r)} w^{-(p'-1)}\right) \frac{dr}{1-r}\leq C
{\rm Re}\, {\mathcal U}_{sp}^{w\lambda}(\nu)(\rho\eta),$$
and letting $\rho\rightarrow 1$, we obtain a).

In order to obtain the norm estimate, lets us simply write ${\mathcal U}(z)= {\mathcal U}_{sp}^{w\lambda} (\nu)(z)$, and prove that for $k>s$,
$$||{\mathcal U}||_{HF_s^{p1}(w)}^p=|{\mathcal U}(0)|^p+\int_{{\bf S}^n} \left( \int_0^1 (1-\rho)^{k-s}|(I+R)^k{\mathcal U}(\rho\eta)|\frac{d\rho}{1-\rho} \right)^pw(\eta)d\sigma(\eta) \leq C{\mathcal E}_{sp}^w(\nu).$$
 
But
\begin{equation*}\begin{split}
& \int_0^1 (1-\rho)^{k-s}|(I+R)^k{\mathcal U}(\rho\eta)|\frac{d\rho}{1-\rho}\\
&\preceq\int_0^1 (1-\rho)^{k-s} \int_0^1 \int_{{\bf S}^n} \left( \frac{\nu(B(\zeta,1-r))}{ (1-r)^{n-sp}}\right)^{p'-1} \times\\&\frac{(1-r)^{\lambda-n}}{|1-\rho r\eta\overline{\zeta}|^{\lambda+k}}\left({\int\!\!\!\!\!\setminus}_{B(\zeta,1-r)} w^{-(p'-1)}\right)d\sigma(\zeta)\frac{dr}{1-r}\frac{d\rho}{ 1-\rho}\preceq\Upsilon(\eta),
\end{split}\end{equation*}
where $$\Upsilon(\eta)= \int_0^1\int_{{\bf S}^n} \left( \frac{\nu(B(\zeta,1-r))}{ (1-r)^{n-sp}}\right)^{p'-1} \frac{(1-r)^{\lambda-n}}{|1- r\eta\overline{\zeta}|^{\lambda+s}}\left({\int\!\!\!\!\!\setminus}_{B(\zeta,1-r)} w^{-(p'-1)}\right)d\sigma(\zeta)\frac{dr}{1-r}.$$
Observe that $|{\mathcal U}(0) |^p \leq C||\Upsilon||_{ L^p(w)}^p$.  
Consequently, in order to finish the proof of the theorem, we just need to show that
\begin{equation}\label{reformulation}
||\Upsilon||_{ L^p(w)}^p\leq C{\mathcal E}_{sp}^w(\nu).
\end{equation}

 H\"older's inequality
  with exponent $\frac1{p-1}>1$ gives that
\begin{equation}\label{holders}
\Upsilon(\eta)\leq \Upsilon_1(\eta)^{p-1}\Upsilon_2(\eta)^{2-p},
\end{equation}
where 
$$
\Upsilon_1(\eta)=\int_0^1 \int_{{\bf S}^n} \frac{\nu(B(\zeta,1-r))}{(1-r)^{n-s}} \frac{(1-r)^{t-n}}{|1-r\eta\overline{\zeta}|^{t} }\left({\int\!\!\!\!\!\setminus}_{B(\zeta,1-r)} w^{-(p'-1)}\right)^{p-1}d\sigma(\zeta) \frac{dr}{1-r},
$$
and
$$\Upsilon_2(\eta)= \int_0^1 \int_{{\bf S}^n}\left( \frac{\nu(B(\zeta,1-r))}{(1-r)^{n-s}}\right)^{p'} \frac{(1-r)^{\frac{\lambda+s-t(p-1)}{2-p}-n}} {|1-r\eta\overline{ \zeta}|^{\frac{\lambda+s-t(p-1)}{2-p}}} \left({\int\!\!\!\!\!\setminus}_{B(\zeta,1-r)} w^{-(p'-1)}\right)^{p} \frac{d\sigma(\zeta)dr}{1-r}.
$$

We begin estimating the function $\Upsilon_1$.
If  $\zeta\in B(\tau,1-r)$, we have that $B(\zeta,1-r)\subset B(\tau, C(1-r))$, and since $w^{-(p'-1)}$ satisfies a doubling condition,
\begin{equation} \label{psi-1}
\Upsilon_1(\eta)\preceq 
\int_0^1 (1-r)^{t-2n+s}\int_{{\bf S}^n}\int_{B(\tau,C(1-r))}     \frac{d\sigma(\zeta)} {|1-r\eta\overline{\zeta}|^{t}}\left({\int\!\!\!\!\!\setminus}_{B(\tau,1-r)} w^{-(p'-1)}\right)^{p-1}\frac{d\nu(\tau) dr}{1-r}.
 \end{equation}
Next, we observe that if $\zeta\in B(\tau,C(1-r))$, $|1-r\eta\overline{\tau}|\preceq |1-r\eta\overline{\zeta}|$. Hence, the above is bounded by
$$
C \int_0^1 (1-r)^{t-n+s}\int_{{\bf S}^n}\frac{\left({\int\!\!\!\!\!\setminus}_{B(\tau,1-r)} w^{-(p'-1)}\right)^{p-1} } {|1-r\eta\overline{\tau}|^{t}}d\nu(\tau) \frac{dr}{1-r}.
$$
Since
\begin{equation*}\begin{split}
&\int_{{\bf S}^n}\frac{\left({\int\!\!\!\!\!\setminus}_{B(\tau,1-r)} w^{-(p'-1)}\right)^{p-1} } {|1-r\eta\overline{\tau}|^{t}}d\nu(\tau)\preceq 
 \int_{{\bf S}^n}\left({\int\!\!\!\!\!\setminus}_{B(\tau,1-r)} w^{-(p'-1)}\right)^{p-1}\int_{|1-r\eta\overline{\tau}| \leq \delta}\frac{d\delta}{\delta^{t+1}}d\nu(\tau), 
\end{split}\end{equation*}
 the above estimate, together with Fubini's theorem and the fact that $t-n+s>\tau-n$ give that $\Upsilon_1(\eta)$ is bounded by
\begin{equation*}\begin{split}
&C\int_0^1\int_{B(\eta,\delta)} \delta^{t-n+s}\left({\int\!\!\!\!\!\setminus}_{ B(\tau,\delta)} w^{-(p'-1)}\right)^{p-1} d\nu(\tau) \frac{d\delta}{\delta^{t+1}}\\
&\preceq\int_0^1\left({\int\!\!\!\!\!\setminus}_{B(\eta,\delta)} w^{-(p'-1)}\right)^{p-1}\frac{\nu(B(\eta,\delta))}{\delta^{n-s}}\frac{d\delta}{\delta},
 \end{split}\end{equation*}
where we have used the fact that if  $\tau\in B(\eta,\delta)$, then $B(\tau,\delta)\subset B(\eta, C\delta)$, for some $C>0$ and that $w^{-(p'-1)}$ satisfies a doubling condition.

Applying H\"older's inequality with exponent $\frac1{(p-1)^2}>1$, we deduce that
\begin{equation}\begin{split}\label{newestimate}&||\Upsilon||_{ L^p(w)}\\ &\preceq\left( \int_{{\bf S}^n}\left( \int_0^1\left({\int\!\!\!\!\!\setminus}_{B(\eta,1-r)} w^{-(p'-1)}\right)^{p-1}\frac{\nu(B(\eta,\delta))}{\delta^{n-s}}\frac{d\delta}{\delta}\right)^{p'}wd\sigma\right)^{(p-1)^2} \left( \int_{{\bf S}^n} \Upsilon_2w\right)^{p(2-p)} .
\end{split}\end{equation}
Theorem \ref{extenssionwolfftheorem} with $q=1$ gives that the first factor on the right is bounded by
$  C{\mathcal E}_{sp}^w(\nu)^{(p-1)^2}$.

Next we deal with the integral involving $\Upsilon_2$. We recall that  $l=\frac{\lambda+s-t(p-1)}{2-p}-n>\tau -n$. Fubini's theorem gives that
\begin{equation*}\begin{split}
&\int_{{\bf S}^n} \Upsilon_2w\\
&=\int_{{\bf S}^n}\int_0^1\left( \frac{\nu(B(\zeta,1-r))}{(1-r)^{n-s}} \right)^{p'}(1-r)^{l}\left({\int\!\!\!\!\!\setminus}_{B(\zeta,1-r)} w^{-(p'-1)}\right)^{p}\int_{{\bf S}^n}\frac{w(\eta)d\sigma(\eta)}{|1-r\eta\overline{\zeta}|^{l+n}}\frac{d\sigma(\zeta)dr}{1-r}.
\end{split}\end{equation*}
But, as before, since $l>\tau-n$,
$$\int_{{\bf S}^n}\frac{w(\eta)d\sigma(\eta)}{|1-r\eta\overline{\zeta}|^{l+n}}\leq \frac{C}{(1-r)^l} {\int\!\!\!\!\!\setminus}_{B(\zeta,1-r)}w.
$$
The above, together with Fubini's theorem gives that
$$\int_{{\bf S}^n} \Upsilon_2w\preceq 
\int_0^1 \int_{{\bf S}^n}{\int\!\!\!\!\!\setminus}_{B(\eta,1-r)}\left( \frac{\nu(B(\zeta,1-r))}{(1-r)^{n-s}}\right)^{p'}\left({\int\!\!\!\!\!\setminus}_{B(\zeta,1-r)} w^{-(p'-1)}\right)^{p}d\sigma(\zeta) w(\eta)\frac{d\sigma(\eta)dr}{1-r}. 
$$
But if $\zeta\in B(\eta,1-r)$, $B(\zeta,1-r)\subset B(\eta,C(1-r))$, for some $C>0$, and in consequence the above is bounded by
$$C\int_{{\bf S}^n} \int_0^1\left( \frac{\nu(B(\eta,C(1-r)))}{(1-r)^{n-s}}\right)^{p'}\left({\int\!\!\!\!\!\setminus}_{B(\eta,1-r)} w^{-(p'-1)}\right)^{p}\frac{dr}{1-r}w(\eta) d\sigma(\eta).
$$
The change of variables $C(1-r)=y-1$ gives that we can estimate the previous expression by
\begin{equation*}\begin{split}
&C\int_{{\bf S}^n} \int_0^1\left( \frac{\nu(B(\eta,(1-y)))}{(1-y)^{n-s}}\right)^{p'} \left({\int\!\!\!\!\!\setminus}_{B(\eta,1-y)} w^{-(p'-1)}\right)^{p}\frac{dy}{1-y}w(\eta) d\sigma(\eta)\\& +\nu({{\bf S}^n})^{p'} \left( \int_{{\bf S}^n}w^{-\frac1{p-1}}\right)^p=I+II.
\end{split}\end{equation*}
Theorem \ref{extenssionwolfftheorem} gives that $II\preceq C {\mathcal E}_{sp}^w(\nu)$, and Theorem \ref{extenssionwolfftheorem} with $q=p'$
gives that $I\leq C{\mathcal E}_{sp}^w(\nu)$. Consequently,
$\int_{{\bf S}^n} \Upsilon_2w\leq C{\mathcal E}_{sp}^w(\nu)$,
and plugging this estimate in (\ref{newestimate}), we deduce that
$$||\Upsilon||_{ L^p(w)}^p\preceq C {\mathcal E}_{sp}^w(\nu)^{(p-1)^2}{\mathcal E}_{sp}^w(\nu)^{p(2-p)}\simeq{\mathcal E}_{sp}^w(\nu).$$

We now sketch the proof of part (2). We   choose $\lambda>0$ such that $\tau-sp<\lambda<1$, and define ${\mathcal V}_{sp}^{w\lambda}(\nu)(z)$ as in (\ref{holompot2}).
 Let us simplify the notation and just write $ {\mathcal V}(z)={\mathcal V}_{sp}^{w\lambda}(\nu)(z)$.
 Let  $\varepsilon\in\R$   such that $\tau <\varepsilon+n<\lambda+sp$. 
 
 The proof of a) is analogous to the one in case $1<p<2$.
 
For the proof of b), 
 let us consider $k>s$. It will be enough to prove the following: \begin{equation}\begin{split}\label{estimatePhi}&||{\mathcal V}||_{HF_s^{p1}(w)}^p\\&=|{\mathcal V}(0)|^p+\int_{{\bf S}^n} \left( \int_0^1 (1-\rho)^{k-s}|(I+R)^k{\mathcal V}(\rho\zeta)|\frac{d\rho}{1-\rho} \right)^pw(\zeta)d\sigma(\zeta) \leq  C{\mathcal E}_{sp}^w(\nu).\end{split}\end{equation}
   Let us begin with the estimate $|{\mathcal V}(0)|^p\preceq {\mathcal E}_{sp}^w(\nu)$. 
  If $p>2$, H\"older's inequality with exponent $\frac1{p'-1}>1$, gives that
  \begin{equation*}\begin{split}
  &|{\mathcal V}(0)| \leq \left(\int_0^1 \int_{{\bf S}^n} (1-r)^{\varepsilon} \left({\int\!\!\!\!\!\setminus}_{B(\zeta, 1-r)}w^{-(p'-1)}\right)^\frac1{p'-1} d\nu(\zeta)\frac{dr}{1-r}\right)^{p'-1}\times\\
  &\left( \int_0^1 \left( (1-r)^{(p'-1)(\lambda+sp-n-\varepsilon)}\right)^\frac1{2-p'} \frac{dr}{1-r} \right)^{2-p'} \preceq  \nu({\bf S}^n)^{p'-1} {\int\!\!\!\!\!\setminus}_{{\bf S}^n }w^{-{(p'-1)}}.
  \end{split}\end{equation*} 
   The case $p=2$ is proved similarly.
Consequently, for any $p\geq 2$, 
  $$|{\mathcal V}(0)|^p\preceq \nu({\bf S}^n)^{p'} \left({\int\!\!\!\!\!\setminus}_{{\bf S}^n }w^{-{(p'-1)}}\right)^p\leq C{\mathcal E}_{sp}^w(\nu),$$
  where the constant $C$ may depend on $w$.
  
  Following with the estimate of $||{\mathcal V}||_{HF_s^{p1}(w)}$, we recall (for example see \cite{cohnverbitsky1}, Proposition 1.4) that if $k>0$,  $0<\lambda<1$, and $z\in\B^n$,
  $$|(I+R)^k\left(\int_{{\bf S}^n} \frac{d\nu(\zeta)}{(1-z\overline{\zeta})^\lambda} \right)^{p'-1}|\leq C \left(\int_{{\bf S}^n} \frac{d\nu(\zeta)}{|1-z\overline{\zeta}|^\lambda} \right)^{p'-2} \int_{{\bf S}^n} \frac{d\nu(\zeta)}{|1-z\overline{\zeta}|^{\lambda+k}}.$$
   
 Plugging this estimate in (\ref{estimatePhi}) and using that $p'-2\leq 0$, we get
 \begin{equation*}\begin{split}
 &|(I+R)^k{\mathcal V}(\rho\eta)| \\& 
\preceq\int_0^1\int_{1-r<\delta,1-\rho<\delta<3}\frac{(1-r)^{(p'-1)(\lambda+sp-n)}\left( \int_{B(\eta,\delta)} 
\left({\int\!\!\!\!\!\setminus}_{B(\zeta,1-r)} w^{-(p'-1)}\right)^\frac1{p'-1}d\nu(\zeta) \right)^{p'-1}}{\delta^ {\lambda+k+1+(p'-2)\lambda}}\frac{d\delta dr}{1-r}.\end{split}\end{equation*}
 Assume first that $p>2$. Fubini's theorem and H\"older's inequality with exponent $\frac1{p'-1}>1$, gives that the above is bounded by
 \begin{equation}\begin{split}\label{estimacio}
 & \int_{1-\rho}^3 \left( \int_{1-r<\delta<3} (1-r)^{\varepsilon} \int_{B(\eta,\delta)}\left({\int\!\!\!\!\!\setminus}_{B(\zeta,1-r)} w^{-(p'-1)}\right)^\frac1{p'-1}d\nu(\zeta) \frac{dr}{1-r}\right)^{p'-1}\times\\&
\left(\int_{1-r<\delta<3} \left( \frac{(1-r)^{(\lambda+sp-n)(p'-1)-\varepsilon(p'-1)}}{\delta^{\lambda+k+1+(p'-2)\lambda}} \right)^\frac1{2-p'} \frac{dr}{1-r} \right)^{2-p'}d\delta. \end{split}\end{equation}
Next, Fubini's theorem and the fact that $\varepsilon>\tau-n$ give that
\begin{equation*}\begin{split}
&\int_{1-r<\delta}(1-r)^\varepsilon \int_{B(\eta,\delta)} 
\left({\int\!\!\!\!\!\setminus}_{B(\zeta,1-r)} w^{-(p'-1)}\right)^\frac1{p'-1}d\nu(\zeta) \frac{dr}{1-r} \\&\preceq \int_{B(\eta,\delta)}\delta^\varepsilon \left({\int\!\!\!\!\!\setminus}_{B(\zeta,\delta)} w^{-(p'-1)}\right)^\frac1{p'-1} d\nu(\zeta).
\end{split}\end{equation*}
 We also have that since $\lambda+sp-n-\varepsilon>0$,
  (\ref{estimacio}) is bounded by
$$\int_{1-\rho}^3 \left(\int_{B(\eta,\delta)} \left({\int\!\!\!\!\!\setminus}_{B(\zeta,\delta)} w^{-(p'-1)}\right)^\frac1{p'-1}d\nu(\zeta) \right)^{p'-1} \frac{d\delta}{\delta^{(n-sp)(p'-1)+k+1}}.$$

For the case $p=2$, we obtain the same estimate, applying directly    condition  (\ref{diniconditionnearzero}) on (\ref{estimacio}).

Integrating with respect to $\rho$, and applying Fubini's theorem we get
\begin{equation*}\begin{split}
& \int_0^1 (1-\rho)^{k-s} |(I+R)^k{\mathcal V}(\rho\eta)|\frac{d\rho}{1-\rho}\\
&\preceq \int_0^3 \left(\int_{B(\eta,\delta)} \left({\int\!\!\!\!\!\setminus}_{B(\zeta,\delta)} w^{-(p'-1)}\right)^\frac1{p'-1}d\nu(\zeta) \right)^{p'-1} \frac{d\delta}{\delta^{(n-s)(p'-1)+1}},
\end{split}\end{equation*}
since $(n-sp)(p'-1)+s=(n-s)(p'-1)$.
If $\tau\in B(\zeta,\delta)$, and $\zeta\in B(\eta,\delta)$, we have that $\tau\in B(\eta, C\delta)$. The fact that $w^{-(p'-1)}$ satisfies a doubling condition, gives that the last integral is bounded by
$${\displaystyle C\int_0^3 \left( \frac{\nu(B(\eta,\delta)}{\delta^{n-s}} \right)^{p'-1} {\int\!\!\!\!\!\setminus}_{B(\eta,\delta)} w^{-(p'-1)} \frac{d\delta}{\delta}.}$$

Applying  Theorem \ref{extenssionwolfftheorem} with exponent $q=p'-1$, we finally obtain that 
$$\int_{{\bf S}^n} \left(\int_0^1 (1-\rho)^{k-s} |(I+R)^k{\mathcal V}(\rho \eta)|\frac{d\rho}{\rho} \right)^p w(\eta)d\sigma(\eta)\preceq
\int_{{\bf S}^n} {\mathcal W}_{sp}^w(\nu)(\zeta) d\nu (\zeta).\qed
$$

We can now state the characterization of the weighted Carleson measures.
\begin{theorem}\label{weightedtraceinequality}
 Let $1< p<+\infty$, $0<n-sp<1$, $w$  an $A_p$-weight, and $\mu$ a finite positive Borel measure on $\B^n$. Assume that $w$ is in $D_\tau$ for some $0\leq \tau-sp<1$. We then have that the following statements are equivalent:
 
 (i) $||K_\alpha(f)||_{L^p(d\mu)}\leq C||f||_{ L^p(w)}$.
 
 (ii) $||f||_{L^p(d\mu)}\leq C|| f||_{H_s^p(w)}$.
\end{theorem}
 {\bf Proof of Theorem \ref{weightedtraceinequality}:}\par
 Let us show first that $(i)\Longrightarrow (ii)$. Theorem \ref{characterizationhardysobolev} gives that   condition (ii) can be rewritten as
 $$||C_s(g)||_{L^p(d\mu)}\leq C||g||_{ L^p(w)}.$$
 This fact together with the estimate $|C_s(f)|\leq CK_s(|f|)$ finishes the proof of the implication.
 
 Assume now that (ii) holds. Since a measure $\mu$ on $\B^n$ satisfies (i) if and only if (see (\ref{equation3.5})) there exists $C>0$ such that for any open set $G\subset {\bf S}^n$,
 ${\displaystyle\mu(T(G))\leq C C_{sp}^w(G),}$ we will check that this estimate holds.
  Let $G\subset {\bf S}^n$ be an open set, and let $\nu$ be the extremal measure for $C_{sp}^w(G)$. We then have that  ${\mathcal W}_{sp}^w(\nu) \geq 1$ except on a set of $C_{sp}^w$-capacity zero, and $\int_{{\bf S}^n}{\mathcal W}_{sp}^w(\nu)d\nu \leq C C_{sp}^w(G)$. Let us check that the first estimate also holds for a.e. $x\in G$ (with respect to Lebesgue measure on ${\bf S}^n$). Indeed, if $A\subset {\bf S}^n$ satisfies that $C_{sp}^w(A)=0$, and $\varepsilon>0$, let $f\geq 0$ be a function such that $K_s(f)\geq 1$ on $A$ and $\int_{{\bf S}^n}f^pw\leq \varepsilon$. Since $ L^p(w)\subset L^{p_1}(d\sigma)$, for some $1<p_1<p$, (see Lemma \ref{weightedspaces}) we then have
 ${\displaystyle ||f||_{L^{p_1}(d\sigma)}\leq C ||f||_{ L^p(w)}\leq C \varepsilon^\frac1{p}.}$
 Thus $C_{sp_1}(A)=0$, and in particular $|A|=0$.
 
  Following with the proof of the implication consider   the holomorphic function on $\B^n$ defined by $F(z)= {\mathcal U}_{sp}^{w\lambda}(\nu)(z)$ if $1<p<2$, $F(z)={\mathcal V}_{sp}^{w\lambda}(\nu)(z)$, if $p\geq 2$ where  $ \lambda$ is  as in  Theorem \ref{pgreaterthan2}. Theorem   \ref{pgreaterthan2} and the fact that $\nu$ is extremal give that
 $$\lim_{r\rightarrow 1} {\rm Re}\, F(r\zeta)\geq C{\mathcal W}_{sp}^w(\nu)(\zeta)\geq C,$$
 for a.e. $x\in G$ with respect to $C_{sp}^w$, and in consequence, for a.e. $x\in G$ with respect to Lebesgue measure on $G$. Hence, if $P$ is the Poisson-Szeg\"o kernel
 $$ |F(z)|=|P[\lim_{r\rightarrow1} F(r\cdot)](z)|\geq |P[{\rm Re}\, \lim_{r\rightarrow 1} F(r\cdot)](z)| \geq C,$$
 for any $z\in T(G)$,  and since we are assuming that (ii) holds, we obtain
 $$\mu(T(E))\leq \int_{T(E)}|F(z)|^pd\mu(z) \leq C|| F||_{H_s^p(w)}^p\leq C {\mathcal E}_{sp}^w (\nu) \leq C C_{sp}^w(G).\qed 
 $$
 
 We finish with an example which shows that, simmilarly to what happens if $w\equiv 1$, if $w\in D_\tau$  and $\tau-sp>1$, then the equivalence between (i) and (ii) in the previous theorem need not to be true.
    
 \begin{proposition}\label{contraexemple}
 Let  $n\geq 1$, $p=2$, and $\tau\geq 0$, $0<s $ such that $\tau>1+2s$. Assume also that $n<\tau<n+1$. Then there exists $w\in A_2\cap D_\tau$ and a positive Borel measure $\mu$ on ${\bf S}^n$ such that  $\mu$ is a Carleson measure for $H_{s}^{2}(w )$, but it is not Carleson for $K_{s }[L^2(w)]$. 
 \end{proposition}
{\bf Proof of Proposition \ref{contraexemple}:}\par

 If $\varepsilon=\tau-n$, and $\zeta=(\zeta',\zeta_n)\in {\bf S}^n$, we consider the weight on ${\bf S}^n$ defined by $w(\zeta)=(1-|\zeta'|^2)^\varepsilon$. A calculation gives that $w(z)=(1-|z|^2)^\varepsilon\in A_2$ if and only if  $-1<\varepsilon<1$, which is our case. We also have that if $\zeta\in{\bf S}^n$, $R>0$ and $j\geq 0$, then $W(B(\zeta,2^jR))\simeq 2^{j\tau}W(B(\zeta,R))$, i.e. $w\in D_\tau$.
 
 Next, any function in $H_s^2(w)$ can be written as $\int_{{\bf S}^n} \frac{f(\zeta)}{(1-z\overline{\zeta})^{n-s}}d\sigma(\zeta)$, $f\in L^2(w)$.  It is then immediate to check that the restriction to $B^{n-1}$ of any such function can be written as
 $$\int_{\B^{n-1}} \frac{g(\zeta')(1-|\zeta'|^2)^{-\frac{\varepsilon}2}}{(1-z'\overline{\zeta'})^{n-s}}dv(\zeta'),$$
 with $g\in L^2(dv)$. 
  This last space coincides (see for instance \cite{peloso}) with the Besov space $B_{s-\frac12-\frac{\varepsilon}2}^2(\B^{n-1})=H_{s-\frac12-\frac{\varepsilon}2}^2({\B}^{n-1})$.
  
Next, $n-1-(s-\frac12-\frac{\varepsilon}2)2=\tau-2s>1$, and Proposition 3.1 in \cite{cascanteortega2} gives that there exists a positive Borel measure $\mu$ on $\B^n$ which is Carleson for $H_{s-\frac12-\frac{\varepsilon}2}^2({\bf S}^{n-1})$, but it fails to be Carleson for the space $ K_{s-\frac12-\frac{\varepsilon}2}[L^2(d\sigma)] $. Thus the operator $$f\longrightarrow \int_{{\bf S}^{n-1}} \frac{f(\zeta)}{|1-z\overline{\zeta}|^{n-1-(s-\frac12-\frac{\varepsilon}2)}}d\sigma(\zeta),$$
is not bounded from $L^2(d\sigma)$ to $L^2(d\mu)$. Duality gives that   the operator
$$g\longrightarrow \int_{\B^{n-1}} \frac{g(z)}{|1-z\overline{\zeta}|^{n-1-(s-\frac12-\frac{\varepsilon}2)}}d\mu(z)$$
is also not bounded from $L^2(d\mu)$ to $L^2(d\sigma)$. But if $g\geq 0$, $g\in L^2(d\mu)$, Fubini's theorem gives
\begin{equation*}\begin{split} 
&||\int_{\B^{n-1}} \frac{g(z)}{|1-z\overline{\zeta}|^{n-1-(s-\frac12-\frac{\varepsilon}2)}}||_{L^2(d\sigma)}^2=
\int_{{\bf S}^{n-1}}\left(\int_{\B^{n-1}} \frac{g(z)}{|1-z\overline{\zeta}|^{n-1-(s-\frac12-\frac{\varepsilon}2)}}d\mu(z)\right)^2d\sigma(\zeta)\\&
=\int_{{\bf S}^{n-1}} \int_{\B^{n-1}} \frac{g(z)}{|1-z\overline{\zeta}|^{n-1-(s-\frac12-\frac{\varepsilon}2)}}d\mu(z)  \int_{\B^{n-1}} \frac{g(w)}{|1-w\overline{\zeta}|^{n-1-(s-\frac12-\frac{\varepsilon}2)}}d\mu(w)d\sigma(\zeta)\simeq\\&
\int_{\B^{n-1}} \int_{\B^{n-1}}\frac{g(z)g(w)}{|1-z\overline{w}|^{n-1-2(s-\frac12-\frac{\varepsilon}2)}}d\mu(z)d\mu(w),
\end{split}\end{equation*}
where   the last estimate holds since $n-1-2(s-\frac12-\frac{\varepsilon}2)=\tau-2s >0$. Consequently, we have that for the measure $\mu$, it does not hold that for any  $g\in L^2(d\mu)$
\begin{equation}\label{convolution}  \int_{\B^{n-1}}\int_{\B^{n-1}}\frac{g(z)g(w)}{|1-z\overline{w}|^{n-2(s-\frac{\varepsilon}2)}}d\mu(z)d\mu(w)\leq C||g||_{L^2(d\mu)}.\end{equation}    

We next check that the failure of being a Carleson measure for $ K_s[L^2(w)] $ can be also rewritten in the same terms. An argument similar to the previous one,  gives that $\mu$ is not Carleson   for $ K_s[L^2(w)] $   if and only if the operator
$$f\longrightarrow \int_{\B^{n-1}} \frac{f(z)}{|1-y\overline{z}|^{n -s}}dv(z)$$
is not bounded from $L^2(wdv)$ to $L^2(d\mu)$. Equivalently, writing $f(z)=h(z)(1-|z|^2)^\frac{\varepsilon}2$,  this last assertion holds if and only if the operator
  $$f\longrightarrow \int_{\B^{n-1}} \frac{f(z)(1-|z|^2)^\frac{-\varepsilon}2}{|1-y\overline{z}|^{n -s}}dv(z)$$
is not bounded from $L^2(dv)$ to $L^2(d\mu)$.
But an argument as before, using duality and Fubini's theorem, gives that the fact that of the unboundedness of the operator can be rewritten in terms of (\ref{convolution}).\qed


\begin{thebibliography}{BZ}

 

\bibitem[{\bf Ad}]{adams} D.R. Adams,   {\it Weighted nonlinear potential
theory}, Trans. Amer. Math. Soc.  {\bf 297} (1986), 73--94.
 
\bibitem[{\bf AdHe}]{adamshedberg} D.R. Adams and L.I. Hedberg, {\it Function
Spaces and Potential Theory}, Springer-Verlag
Berlin--Heidelberg--New York, 1996.
\bibitem[{\bf Ah}]{ahern} P. Ahern, {\it Exceptional sets for holomorphic Sobolev functions}, Michigan Math. J.
  {\bf 35 }, (1988),  29--41.
 
\bibitem[{\bf AhCo}]{aherncohn} P. Ahern and W.S. Cohn, {\it Exceptional sets for Hardy-Sobolev spaces  }, Indiana Math. J.
  {\bf 39 }, (1989),  417--451.
  
  \bibitem[{\bf AhBrCa}]{ahernbrunacascante} P. Ahern, J. Bruna and C. Cascante, {\it $H^p$-theory for generalized $M$-harmonic functions in the unit ball }, Indiana Math. J.
  {\bf 45 }, (1996),  103--135.
   
   
 \bibitem[{\bf BeLo}]{berglofstrom} J. Berg and J. L\"ofstr\"om, {\it Interpolation Spaces, an Introduction}, Springer-Verlag
Berlin, 1976.
   

\bibitem[{\bf CaOr1}]{cascanteortegat}
C. Cascante and J.M. Ortega, {\it Tangential-exceptional sets for Hardy-Sobolev spaces}, Illinois
J. Math. {\bf 39}, (1995), 68--85.

\bibitem[{\bf CaOr2}]{cascanteortega2}
C. Cascante and J.M. Ortega, {\it Carleson measures on spaces of Hardy-Sobolev type}, Canadian
J. Math. {\bf 47}, (1995), 1177--1200.
 


 
  
 
\bibitem[{\bf CohVe1}]{cohnverbitsky2} W.S. Cohn, I.E. Verbitsky, {\it  Trace inequalities for Hardy-Sobolev functions in the unit ball of $\C^n$}, Indiana Univ. Math. J.
  {\bf 43 }, (1994),
 1079--1097.
 \bibitem[{\bf CohVe2}]{cohnverbitsky1} W.S. Cohn, I.E. Verbitsky, {\it  Non-linear
 potential theory on the ball, with applications to exceptional and boundary interpolation sets}, Michigan Math. J.
  {\bf 42 }, (1995),
 79--97.
\bibitem[{\bf CoiMeSt}]{coifmanmeyerstein} R.R. Coifman, Y. Meyer and E.M. Stein, {\it  Some new function spaces and their applications to harmonic analysis}, Journal of Funct. Anal.
  {\bf 62 }, (1985),
 304--335.



 
 \bibitem[{\bf HeWo}]{hedbergwolff}
L.I. Hedberg and Th. H. Wolff, {\it Thin sets in nonlinear
potential theory},
 Ann. Inst. Fourier (Grenoble)
{\bf 33}, (1983), 161--187.
   
 

\bibitem[{\bf KaKo}]{kangkoo}
H. Kang and H. Koo, {\it Two-weighted inequalities for the derivatives of holomorphic functions and Carleson measures on the ball},
 Nagoya Math. J.,
{\bf 158}, (2000), 107--131.
\bibitem[{\bf KeSa}]{kermansawyer}
R. Kerman and E.T. Sawyer, {\it  The trace inequality and eigenvalue estimates for Schr\"odinger operators},
 Ann. Inst. Fourier,
{\bf 36}, (1986), 207--228.

 
\bibitem[{\bf Lu}]{luecking}
D. H. Luecking, {\it Representation and duality in weighted spaces of analytic functions},
 Indiana Univ. Math.  
{\bf 34}, (1985), 319--336.
\bibitem[{\bf Ma}]{mazya} V. G. Maz'ya, {\it Sobolev Spaces},
Berlin: Springer, 1985.
\bibitem[{\bf OF}]{ortegafabrega}
J.M. Ortega   and J. Fabrega, {\it Holomorphic Triebel-Lizorkin Spaces},
J. Funct. Analysis
 {\bf 151}, (1997), 177--212.
 
 \bibitem[{\bf Pe}]{peloso}
M. M. Peloso, {\it M\"obius invariant spaces on the unit ball},
Michigan Math. J. {\bf 39} (1992), 509--537.

\bibitem[{\bf Ru}]{rudin} W. Rudin, {\it Function Theory in the Unit Ball of $\C^n$},
New York: Springer, 1980.

 


 
\bibitem[{\bf St2}]{stein2} E.M.Stein, {\it Boundary behavior of holomorphic functions of several complex variables},
Princeton University Press, 1972.
\bibitem[{\bf StrTo}]{strombergtorchinsky} J.-O. Str\"omberg and A. Torchinsky, {\it Weighted Hardy Spaces},
Lecture Notes in Math. 1381, Springer-Verlag 1989.
 

 



\end{thebibliography}
\end{document}